\documentclass{amsart}

\usepackage{amsmath, amsthm}
\usepackage{amssymb}
\usepackage{amsfonts}
\usepackage{latexsym}
\usepackage{mathrsfs}
\usepackage{bm}
\usepackage{graphicx}
\usepackage{enumerate}
\usepackage[utf8]{inputenc}
\usepackage{mathtools}
\usepackage{esint}

\theoremstyle{plain}
\newtheorem{theorem}{Theorem}[section]
\newtheorem{cor}[theorem]{Corollary}
\newtheorem{lem}[theorem]{Lemma}
\newtheorem{pro}[theorem]{Proposition}
\newtheorem{Def}[theorem]{Definition}
\newtheorem{rem}[theorem]{Remark}
\numberwithin{equation}{section}

\newtheorem{thx}{Theorem}

\newcommand{\aas}{along a subsequence}

\newcommand{\hym}{hyperbolic metric}
\newcommand{\ub}{uniformly bounded}
\newcommand{\afm}{almost Fuchsian manifold}

\newcommand{\wrt}{with respect to}
\newcommand{\TS}{Teichm\"{u}ller space}

\newcommand{\cp}{critical point}
\newcommand{\ifif}{if and only if}

\newcommand{\maxp}{maximum principle}
\newcommand{\mtps}{mountain pass}
\newcommand{\mpsl}{mountain pass solution}
\newcommand{\hqd}{holomorphic quadratic differential}
\newcommand{\RS}{Riemann surface}

\newcommand{\sff}{second fundamental form}
\newcommand{\gse}{Gauss equation}
\newcommand{\af}{almost Fuchsian}

\newcommand{\mt}{Moser-Trudinger inequality}

\newcommand{\htm}{hyperbolic three-manifold}
\newcommand{\ms}{minimal surface}
\newcommand{\mi}{minimal immersion}
\newcommand{\cs}{conformal structure}

\newcommand{\R}{\mathbb{R}}

\newcommand{\Ical}{\mathcal{I}}
\newcommand{\Jcal}{\mathcal{J}}
\newcommand{\Ccal}{\mathcal{C}}
\newcommand{\ttc}{\mathtt{c}}
\newcommand{\Zcal}{\mathcal{Z}}
\newcommand{\Pcal}{\mathcal{P}_t}

\newcommand{\be}{\begin{equation}}
\newcommand{\ene}{\end{equation}}
\newcommand{\br}{\begin{rem}}
\newcommand{\er}{\end{rem}}
\newcommand{\bl}{\begin{lem}}
\newcommand{\bcor}{\begin{cor}}
\newcommand{\ecor}{\end{cor}}
\newcommand{\el}{\end{lem}}
\newcommand{\bd}{\begin{Def}}
\newcommand{\ed}{\end{Def}}
\newcommand{\ben}{\begin{enumerate}}
\newcommand{\een}{\end{enumerate}}
\newcommand{\bp}{\begin{proof}}
\newcommand{\ep}{\end{proof}}
\newcommand{\bpo}{\begin{pro}}
\newcommand{\epo}{\end{pro}}
\newcommand{\beq}{\begin{equation*}}
\newcommand{\eeq}{\end{equation*}}
\newcommand{\bea}{\begin{eqnarray}}
\newcommand{\eea}{\end{eqnarray}}
\newcommand{\bear}{\begin{eqnarray*}}
\newcommand{\eear}{\end{eqnarray*}}
\newcommand{\bt}{\begin{theorem}}
\newcommand{\et}{\end{theorem}}

\numberwithin{equation}{section}

\allowdisplaybreaks

\def\XXint#1#2#3{{\setbox0=\hbox{$#1{#2#3}{\int}$}
    \vcenter{\hbox{$#2#3$}}\kern-.5\wd0}}

\makeatletter
\def\@citestyle{\m@th\upshape\mdseries}
\def\citeform#1{{\bfseries#1}}
\def\@cite#1#2{{%
  \@citestyle[\citeform{#1}\if@tempswa, #2\fi]}}
\@ifundefined{cite }{%
  \expandafter\let\csname cite \endcsname\cite
  \edef\cite{\@nx\protect\@xp\@nx\csname cite \endcsname}%
}{}
\makeatother
\begin{document}

\title
{Bifurcation for Minimal Surface Equation in Hyperbolic $3$-Manifolds}

\author{Zheng Huang}
\address{The Graduate Center, The City University of New York,
365 Fifth Ave., New York, NY 10016, USA}
\address{Department of Mathematics,
The City University of New York,
Staten Island, NY 10314, USA.}
\email{zheng.huang@csi.cuny.edu}

\author{Marcello Lucia}
\address{Department of Mathematics,
The City University of New York,
Staten Island, NY 10314, USA.}
\email{marcello.lucia@csi.cuny.edu}

\author{Gabriella Tarantello}
\address{Dipartimento di Matematica, Universita' di Roma ``Tor Vergata", Via della Ricerca Scientifica, I-00133 Roma, ITALY.}
\email{tarantel@mat.uniroma2.it}

\date{May 02, 2019}

\subjclass[2000]{Primary 53C21, Secondary 35J20, 53A10}

\begin{abstract}
Initiated by the work of Uhlenbeck in late 1970s, we study existence, multiplicity and asymptotic behavior for {\mi}s of a closed 
surface in some {\htm}, with prescribed conformal structure on the surface and {\sff} of the immersion. We prove several results in these directions, 
by analyzing the Gauss equation governing the immersion. We determine when existence holds, and obtain unique stable solutions for area minimizing 
immersions. Furthermore, we find exactly when other (unstable) solutions exist and study how they blow-up. We prove our class of unstable solutions 
exhibit different blow-up behaviors when the surface is of genus two or greater. We establish similar results for the blow-up behavior of any general 
family of unstable solutions. This information allows us to consider similar {\mi} problems when the total extrinsic curvature is also prescribed. 
\end{abstract}

\maketitle

\setcounter{section}{-1}

\section {Introduction: Geometric Settings}

Minimal surfaces have long been a fundamental object of intense study in geometry and analysis. In this paper we study {\mi}s of a closed surface
in some {\htm}s.  Inspired by Uhlenbeck's approach (\cite{Uhl83}), results on existence and multiplicity of such {\mi}s, as well as their geometrical 
interpretations, are obtained by analyzing bifurcation properties of solutions to the {\ms} equation. Throughout the paper, we assume $S$ is a closed oriented 
surface of genus $g \geq 2$. The {\TS} of $S$ is denoted by $T_g(S)$, and it is the space of {\cs}s (or equivalently {\hym}s) on $S$ such that 
two {\cs}s are equivalent if there is between them an orientation-preserving diffeomorphism in the homotopy class of the identity. 

%The existence and multiplicity questions of closed {\ms}s in {\htm}s have had important applications in many related fields such as three-manifold topology, 
%Teichm\"uller theory, representation theory of surface groups in $PSL(2,\mathbb{C})$ (see for instance \cite{FHS83, Rub05} and many others).

 When $S$ is immersed in some {\htm} $M$, we denote by $g_0$ the induced metric from the immersion and by $\sigma \in T_g(S)$ the {\cs} on $S$ 
 induced by $g_0$. Furthermore, we denote by $g_\sigma$ the unique {\hym} on $(S,\sigma)$, and by $dA$ its area form. Since the metrics $g_\sigma$ and 
 $g_0$ are conformally equivalent, for a suitable conformal factor $u\in C^\infty(S)$, we have 
\be\label{confmetric}
g_0 = e^{2u}g_\sigma.
\ene 

We denote always by $z = x+iy$ the conformal coordinates on $(S,\sigma)$. So in local conformal coordinates we may write:
\beq
g_\sigma = e^{2u_\sigma}dzd\bar{z}, \ \ \ \text{ and } \ \ \ \ g_0 = e^{2v}dzd\bar{z},
\eeq 
with $v = u_\sigma + u$ and $u$ given in \eqref{confmetric}. Now, in such coordinates, the {\sff} II takes the following quadratic expression:
\be\label{sff0}
\text{II} = h_{11}(dx)^2+ 2h_{12}dxdy+ h_{22}(dy)^2,
\ene
with $h_{11}=-h_{22}$ accounting for the fact that $(S, g_0)$ is a {\ms} in $M$.

The Riemann curvature tensor $R_{ij\kappa\ell}$ and the metric tensor $g = (g_{ij})$ of the {\htm} $(M, g)$ satisfy the following equations: 
\be\label{RieTen}
R_{ij\kappa\ell} = -(g_{i\kappa}g_{j\ell} - g_{i\ell}g_{j\kappa}).
\ene
Note that by Bianchi identities, only six components of $R_{ij\kappa\ell}$ are independent. 
%And we require the induced metric on $S$ from the {\mi} lies in the conformal structure $\sigma$.

In this respect, we can use normal coordinates $(z,r)$ for the normal bundle $T^N(S)$, with conformal coordinates $z \in S$ and 
$r \in (-a,a)$ for some $a> 0$ small. We obtain via the exponential map a local coordinate system on $M$ around $S$, where we have 
$g_{j3} = \delta_{j3}$, for $j = 1,2,3$. In these coordinates, the remaining components $g_{i\kappa}$, $1 \le i\le \kappa \le 2$, are just 
$-R_{i3\kappa3}$ in view of \eqref{RieTen}, namely,
\be\label{RieTen3}
R_{i3\kappa3} = -g_{i\kappa}.
\ene
The equations in \eqref{RieTen3} can be viewed as a second order system of ODEs for $g_{ik}$ in the variable $r$ (and fixed $z \in S$). So 
we can uniquely identify $g_{i\kappa}$ (around $S$) by its initial data: 
\be\label{0-4a}\left\{
  \begin{array}{cl}
    g_{i\kappa}(z,0) = (g_{0})_{i\kappa}(z)
    \vspace{5pt} \\
    \frac12\frac{\partial}{\partial r}g_{i\kappa}(z,0) = h_{i\kappa}(z), \ \ \ 1\le i \le \kappa \le 2.                
  \end{array} \right.
  \qquad
\ene

Such initial data on $S$ are provided by the remaining equations in \eqref{RieTen}. To verify this, we take $\ell =3$ 
and $j\not=3$ in \eqref{RieTen} and get 
\be\label{0-5a}
R_{ij\kappa3} = 0,
\ene
which expresses the Codazzi equations on $S$, for $1\le i, j, \kappa \le 2$. Again only two of those equations are independent, and they 
ensure exactly that the quadratic differential $\alpha = (h_{11} - ih_{12})dz^2$ is holomorphic and
\be\label{sff1}
\text{II} = Re(\alpha),
\ene
see [\cite{Hop89,LJ70}]. In other words, $\alpha \in Q(\sigma)$, where we denote by $Q(\sigma)$ the space of {\hqd}s on $(S,\sigma)$.

Finally taking $i=\kappa=1$ and $j=\ell =2$ in \eqref{RieTen}, one gets: 
\be\label{0-5b}
R_{1212} = -g_{11}g_{22} + g_{12}^2,
\ene
which simply gives the {\gse} on $S$, and it states that the conformal factor $u(z)$ in \eqref{confmetric} on $S$ must satisfy:
\be\label{ge}
\Delta u + 1 -e^{2u} - \frac{|\alpha|^2}{det(g_\sigma)}e^{-2u} = 0,
\ene
with $\alpha$ given in \eqref{sff1}, and $\Delta$ the Laplacian in the {\hym} $g_\sigma$. Indeed, \eqref{ge} simply expresses a consistency 
condition on $(S,g_\sigma)$ between the intrinsic curvature of the metric $g_0$ and the extrinsic curvature $det_{g_0}$II, see \cite{Uhl83} 
for details.

Note that, by Bianchi identities, once the equations \eqref{0-5a}, \eqref{0-5b} hold on $S$ then they hold throughout the normal bundle of $S$. 
Furthermore, these equations on $S$ provide the initial data \eqref{0-4a} in terms of $(\sigma, \alpha)$, simply by using the solutions of the 
Codazzi-Gauss equations \eqref{sff1}, \eqref{ge} into  \eqref{confmetric}-\eqref{sff0}.

Thus by prescribing $\sigma\in T_g(S)$ and $\alpha\in Q(\sigma)$ such that the Codazzi-Gauss equations \eqref{sff1}, \eqref{ge} are solvable, 
it is natural to ask whether it is possible to obtain a {\mi} of $(S,\sigma)$ into some {\htm} with the {\sff} satisfying \eqref{sff1}. In short, we shall 
call a {\it {\mi} of $S$ with data $(\sigma,\alpha)$} any of such {\mi}.

A general construction of a {\mi} with prescribed data satisfying the Codazzi-Gauss equations (called ``hyperbolic germs" in \cite{Tau04}) is 
available in literature, see for instance  \cite{Tau04, Jac82}. However, it is not always possible to guarantee that the corresponding {\htm} is 
complete, unless we are more specific about the induced metric $g_0$ or equivalently about the solution of \eqref{ge}. Thus, to obtain more 
satisfactory results of geometrical nature, Uhlenbeck analyzed in \cite{Uhl83} more closely the set of solutions of \eqref{ge}.

We recall that a solution $u$ of \eqref{ge} is called {\it stable} if the linearized operator of \eqref{ge} at $u$ is nonnegative definite in $H^1(S)$, and 
called {\it strictly stable} if the linearized operator of \eqref{ge} at $u$ is positive definite in $H^1(S)$. The interest to stable solutions is justified by the 
fact that they give rise to (local) area minimizing immersions.

By setting: 
\beq
\|\alpha\|_\sigma^2 =\frac{|\alpha|^2}{det(g_\sigma)} =\frac12\|\text{II}\|_{g_\sigma}^2,
\eeq
the length squared of $\alpha$ {\wrt} the {\hym} $g_\sigma$, Uhlenbeck (\cite{Uhl83}) considered a one-parameter family of {\gse}s:
\be\label{t-ge}
\Delta u + 1 -e^{2u} - t^2\|\alpha\|_\sigma^2e^{-2u} = 0,
\ene
for {\mi}s of $S$ with data $(\sigma,t\alpha)$. Using the implicit function theorem, she proved the existence and uniqueness of a smooth solution curve 
of stable solutions to the equation \eqref{t-ge}:
\vskip 0.1in
\noindent
\bt\label{Uh}(\cite {Uhl83}) Fixing a {\cs} $\sigma \in T_g(S)$, and $\alpha \in Q(\sigma)$, there exists a constant $\tau_0 > 0$, depending only on 
$(\sigma, \alpha)$, such that the equation \eqref{t-ge} admits a unique stable solution if and only if $t\in [0,\tau_0]$. 
Furthermore for each $t\in [0,\tau_0]$, the stable solution $u_t<0$ of \eqref{t-ge} forms a smooth monotone decreasing curve {\wrt} $t$. 
Moreover, $u_t$ is strictly stable for $t\in [0,\tau_0)$ and $u_t \nearrow u_{t=0}=0$, as $t\searrow 0$, in $H^1(S)$.
\et
\vskip 0.1in
From Theorem ~\ref{Uh}, a bifurcation diagram (especially for the lower branch, in absolute value, of the solution curve) can be sketched as below in Figure 1:
%\newpage
\begin{figure}[htbp]
\begin{center}
   \includegraphics[scale=1.3]{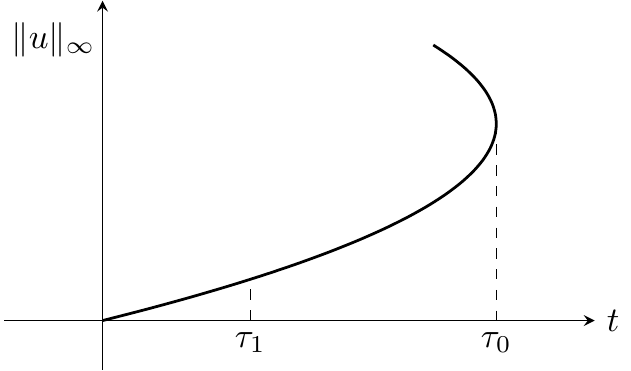}
\end{center}
\caption{Uhlenbeck's Solution Curve}
\end{figure}
\vskip 0.2in
In this diagram, Uhlenbeck indicated a first turn of the curve of stable solutions at some $\tau_0$, though it is still possible the curve retracts and passes 
again the value $t = \tau_0$ to join other solutions of \eqref{t-ge} for $t\ge \tau_0$. Actually it is our first task here to show that this is never the case.

From the geometrical point of view, by Theorem ~\ref{Uh} we know that, for $t\in [0,\tau_0]$, there exists an area minimizing (stable) immersion of 
$S$ with data $(\sigma, t\alpha)$ whose induced metric on the surface $(S,\sigma)$ is $g_0 = e^{2u_t}g_\sigma$. The second variation of the area 
functional is explicitly computed in terms of the linearized operator in \cite {Uhl83}. Also it is interesting to note that there exists a $\tau_1>0$, such that 
for each $t\in (0,\tau_1)$, the surface can be minimally immersed into a so-called {\afm} and this {\af} manifold contains $(S,\sigma)$ as its unique {\ms}. 
In particular, as $t\to 0^{+}$ such family of {\af} {\htm}s converge (in the sense of Gromov-Hausdorff) to a Fuchsian manifold where $S$ is embedded as a 
totally geodesic area minimizing surface, see \cite{Uhl83} for details.

Further work in \cite{HL12} obtained an additional solution for each of Uhlenbeck's (strictly) stable solution to the {\gse}:
\vskip 0.1in
\bt\label{HL}(\cite{HL12}) Let $S$ be a closed surface and $\sigma \in T_g(S)$ be a {\cs} on $S$. If
$\alpha \in Q(\sigma)$ is a {\hqd} on $(S,\sigma)$, then:
\ben[{\bf i)}]
\item
for sufficiently large $t$, the {\gse} \eqref{t-ge} admits no solutions, i.e., there is no {\mi} of $S$ with data
$(\sigma,t\alpha)$ into any {\htm};
\item
for each $t \in (0,\tau_0)$, with $\tau_0>0$ given in Theorem ~\ref{Uh}, there exists also an unstable immersion of $S$ with data $(\sigma,t\alpha)$.
\end{enumerate}
\et
\vskip 0.1in
These results reveal further details on the solution curve to \eqref{t-ge} and an improved bifurcation diagram can be sketched as follows:
\begin{figure}[htbp]
\begin{center}
   \includegraphics[scale=0.5]{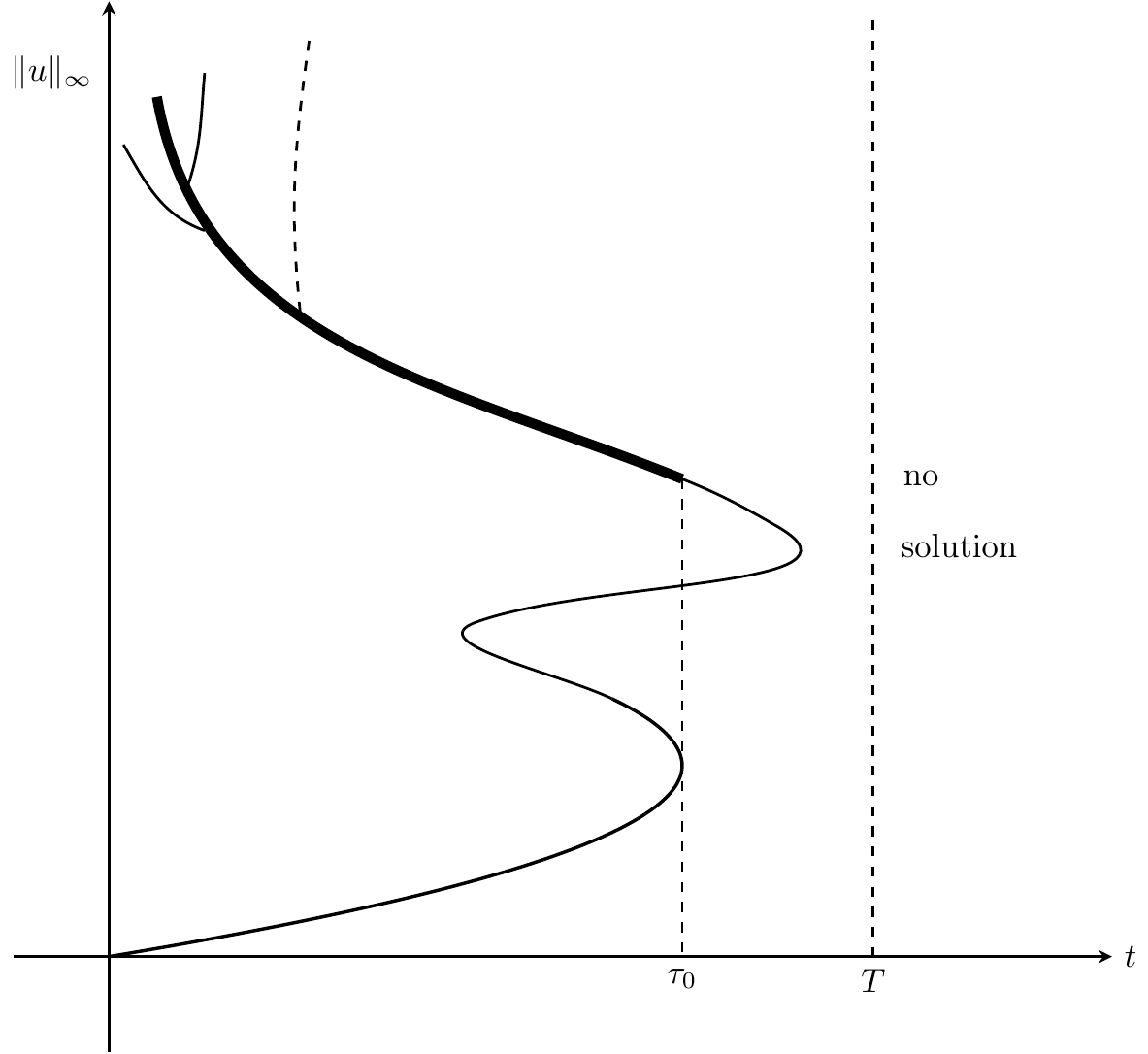}
\end{center}
\caption{Solution Curve from \cite{HL12}}
\end{figure}
%%%%%%%%%%%%%%%%%%%%%%%%%%%%%%%%%%%%%%%%%%%
\section{Introduction: Main results}
%%%%%%%%%%%%%%%%%%%%%%%%%%%%%%%%%%%%
The first purpose of this paper is to complete the above results in Theorems ~\ref{Uh} and ~\ref{HL} as follows. Firstly, we already mentioned, we show 
that actually the interval $[0,\tau_0]$ exhausts the full range of values $t\ge 0$ for which the equation \eqref{t-ge} is solvable. Namely, the bifurcation 
curve starting from the trivial solution at $t=0$, cannot admit an ``S-shape" (as typical in similar problems), but it turns only at $\tau_0$ and until then, 
furnishes the lower branch (in absolute value) of unique (strictly) stable solutions of \eqref{t-ge}. In fact, equation \eqref{t-ge} admits a unique (stable 
but not strictly stable) solution for $t = \tau_0$ and no solutions for $t>\tau_0$. Furthermore we provide a family of unstable solutions for \eqref{t-ge} 
with a specific asymptotic behavior, as $t\to 0^+$. A sketch of the bifurcation diagram can be seen as below in Figure 3.
\begin{figure}[htbp]
\begin{center}
   \includegraphics[scale=0.6]{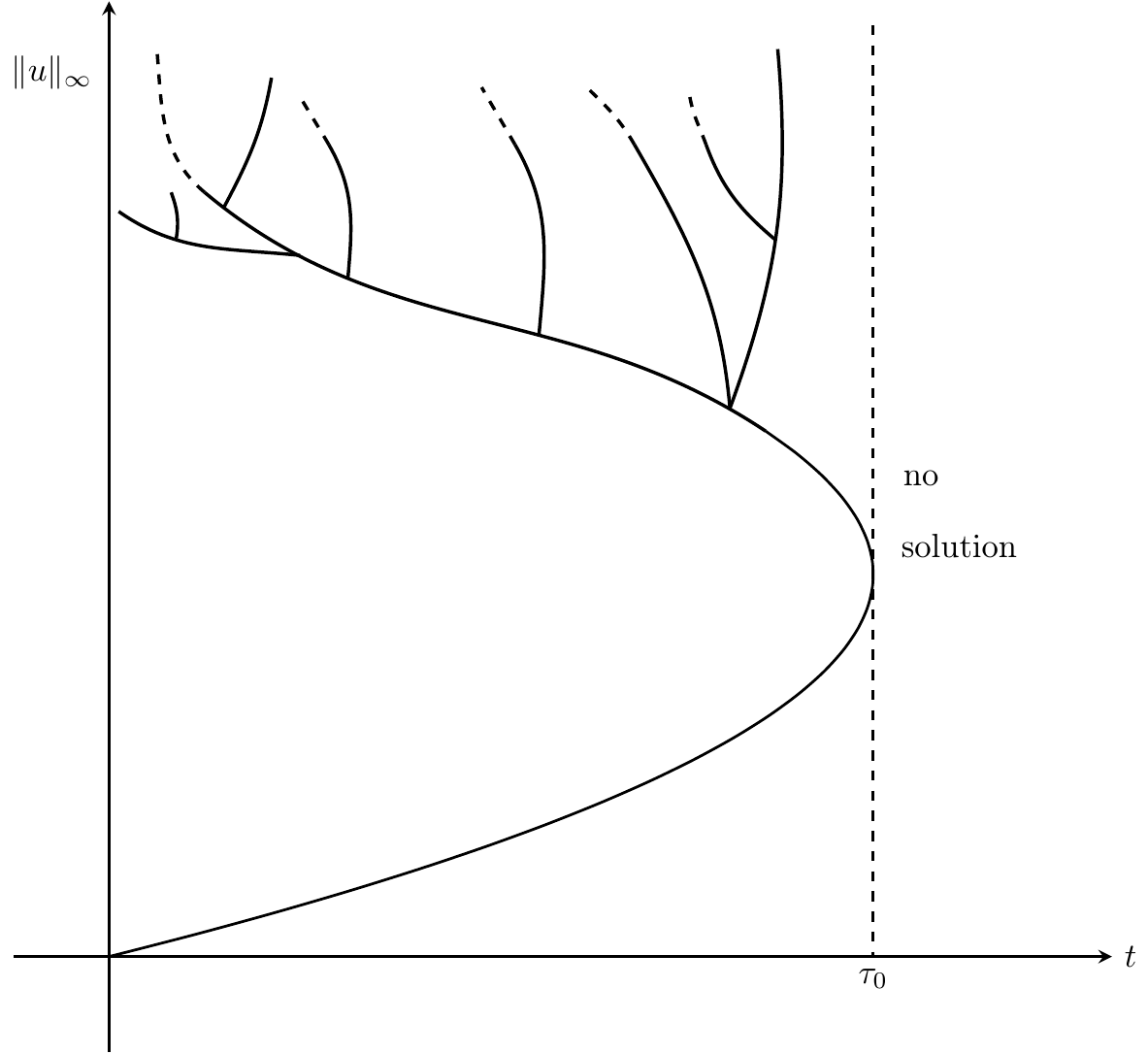}
\end{center}
\caption{New Solution Curve}
\end{figure}
\vskip 0.1in
\begin{thx}\label{main1}
Fixing a {\cs} $\sigma \in T_g(S)$, and a {\hqd} $\alpha \in Q(\sigma)$, the equation \eqref{t-ge} admits a solution if and only if $t \in [0,\tau_0]$, with 
$\tau_0 = \tau_0(\sigma,\alpha)> 0$ given in Theorem ~\ref{Uh}. Furthermore, the unique stable solution $u_t<0$ is the (pointwise) largest solutions of 
\eqref{t-ge} for $t \in [0,\tau_0]$. Moreover, 
\ben[{\bf (i)}]
\item
$\forall t \in (0,\tau_0)$, the equation \eqref{t-ge} admits an \underline{unstable} solution 
$\tilde{u}_t$ (with $\tilde{u}_t <u_t < 0$ on $S$) such that, as $t\searrow 0$, 
\beq
\max\limits_{S}|\tilde{u}_t|\to +\infty,
\eeq
%and 
%\beq
 %t^2\int_S K(z)e^{-2\tilde{u}_t}\ dA \to 4\pi.
 %\eeq
\item
for $t = \tau_0$, the equation \eqref{t-ge} admits the unique solution $u_0$: 
\beq
u_0(z) =\lim_{t\nearrow\tau_0}u_t(z) =\inf\limits_{t\in (0,\tau_0)}u_t(z), \ \forall z\in S
\eeq
\een 
\end{thx}
\vskip 0.1in
We may say that $(\tau_0,u_{\tau_0})$ is a ``bending point" for the bifurcation curve starting at $(t=0,u=0)$, using the terminology introduced in \cite{AR73}. 
Actually we provide a much more detailed study of the asymptotic behavior of the unstable solution $\tilde{u}_t$ in Theorem A, as $t\searrow 0$. Interestingly, 
its behavior depends on whether the surface $(S,\sigma)$ can be ``uniformized" within the class of {\hym}s with conical singularities along a divisor of modulus $2$.

To be more precise, recall that a divisor $D$ on $S$ is a ``formal" expression of the type 
\be\label{divisor}
D = \sum_{j=1}^m\alpha_jp_j,
\ene
with given points $p_j \in S$ and values $\alpha_j > -1, \ j=1,\cdots, m$; and 
\beq
|D| = \sum_{j=1}^m\alpha_j
\eeq
is the modulus of the divisor. Therefore a {\hym} on $S$ with conical singularities along the divisor $D$ in \eqref{divisor} is a metric with Gauss curvature $-1$ 
in $S\backslash \{p_1,\cdots,p_m\}$ and a conical singularity at $p_j$ with angle $\theta_j = 2\pi(1+\alpha_j)$ for $j = 1, \cdots, m$. Thus, the ``uniformization" 
of $(S,\sigma)$ in the class of {\hym}s with conical singularities along the divisor $D$, just means that at least one of such metric is conformal to $g_\sigma$. 
Consequently, the corresponding conformal factor $u$ is given by the \underline{unique} solution of the following equation:
\be\label{conic}
-\Delta u = 1-e^{2u} - 2\pi\sum_{j=1}^m\alpha_j\delta_{p_j}.
\ene   
By analyzing \eqref{conic}, one easily see that the ``uniformization" of $S$ in the sense above is possible if and only if 
\beq
\chi(S) + |D| < 0,
\eeq
where $\chi(S) = 2(1-g)$ is the Euler characteristic of $S$. In particular, if $|D| =2$, then necessarily the genus $g \ge 3$. Thus our blow-up analysis must 
distinguish between the cases where $S$ has genus two or higher.

We start with the following:
\vskip 0.1in
\begin{thx}\label{main2}
Let the genus of the surface $S$ satisfy $g \ge 3$, and $\tilde{u}_t$ be the unstable solution given by Theorem ~\ref{main1}. Then as $t \searrow 0$, we 
have: 
\be\label{blowup3}
t^2\|\alpha\|_\sigma^2e^{-2\tilde{u}_t}\rightharpoonup 4\pi\delta_{p_0}
\ene
weakly in the sense of measures, with some point $p_0 \in S$ such that $\alpha(p_0)\not=0$. Furthermore:
\beq
\tilde{u}_t \to \tilde{u} \ \text{ in } W^{1,q}(S), 1<q<2, \ \text{ uniformly in }C_{loc}^{2,\beta}(S\backslash\{p_0\}), 0<\beta<1,
\eeq
and
\beq
e^{2\tilde{u}_t}\to e^{2\tilde{u}}  \ \ \ \text{ in } L^{s}(S), \ \ \forall s\ge 1,
\eeq
with $\tilde{u}$ the unique solution of the following problem on $S$.:
\be\label{asypsol1}
\Delta\tilde{u} + 1 -e^{2\tilde{u}}-4\pi\delta_{p_0} = 0.
\ene
\end{thx} 
\vskip 0.1in
From the geometrical point of view, Theorem B states that, $\forall t \in (0,\tau_0)$ there exist (unstable) {\mi}s of $S$ with data $(\sigma,t\alpha)$ which 
converge (in the sense of Gromov-Hausdorff) as $t\to 0^{+}$ to a ``limiting" totally geodesic immersion of $S$ into a three-dimensional hyperbolic 
cone-manifold of the type introduced by Krasnov-Schlenker \cite{KS07}. More precisely, the given {\htm} contains conical singularities along one line to infinity 
and the induced metric on $S$ is hyperbolic with a conical singularity along the divisor $D=2p_0$, with suitable $p_0\in S$ and $\alpha(p_0)\not=0$. 
Therefore for our unstable immersions, we observe a quite different limiting behavior from the case of area minimizing (stable) immersions.

If $S$ has genus $g=2$, then the asymptotic behavior of $\tilde{u}_t$ is described by a ``concentration-compactness" alternative, resolved by the existence 
(or not) of a minimum for a Moser-Trudinger type functional. To be more precise, let us recall the following:
\bd
Let $E = \{w \in H^1(S)$ with $\int_Sw(z)dA =0\}$, the Moser-Trudinger functional on $(S,\sigma)$ with weight function $0 \le K \in L^\infty(S)$ is given by   
\be\label{J}
\mathcal{J}(w) = \frac12\int_S|\nabla w|^2\ dA-8\pi\log(\fint_SK(z)e^wdA), \ \ \ w\in E,
\ene
%with weight $K(z) = \|\alpha\|_\sigma^2= \frac{|\alpha|^2}{g_\sigma^2}$, 
where we have used the standard notation:
\beq
\fint_S f\ dA = \frac{\int_S f\ dA}{|S|},
\eeq
and the area $|S| =4\pi(g-1)$ by the Gauss-Bonnet theorem.
\ed
By the Moser-Trudinger inequality (see \cite{Aub98}), the functional $\Jcal$ is bounded from below but \underline{not} coercive in $E$. In other 
words, it is well defined
\be\label{1-1a}
\inf_E\mathcal{J} > -\infty,
\ene
but the infimum in \eqref{1-1a} may not be attained.

For the case of genus $g =2$, our main result reads:
\vskip 0.1in
\begin{thx}\label{main3}
If the surface $S$ is of genus $g=2$, and $\tilde{u}_t$ is the unstable solution given by Theorem ~\ref{main1}, then, as $t\searrow 0$, 
\beq 
\int_S|\tilde{u}_t|\ dA \to +\infty.
\eeq
Furthermore, for $K(z) = \|\alpha\|_\sigma^2$, we have the following alternatives:
\ben[{\bf (i)}]
\item (Compactness) \underline{either}, the Moser-Trudinger functional $\Jcal$ (with $K=\|\alpha\|_\sigma^2$) attains its infimum in $E$, 
and along a sequence $t_n\to 0$, there holds:
\beq
(\tilde{u}_{t_n} - \fint_S \tilde{u}_{t_n}) \to \hat{w}_0, \ \ \ \text{ strongly in }\ H^1(S),
\eeq
and,
\beq
t_n^2K(z)e^{-2\tilde{u}_{t_n}} \to \frac{4\pi K(z)e^{-2\hat{w}_0}}{\int_S K(z)e^{-2\hat{w}_0}\ dA} \ \ \ \text{ uniformly in }\ C^{2,\beta}(S),
\eeq
with $\hat{w}_0$ satisfying on $(S,\sigma)$,  
\be\label{asypsol2}
  \left\{
  \begin{array}{cl}
  \Delta \hat{w}_0 + 4\pi \left(\frac{1}{|S|} - \frac{K(z)e^{-2\hat{w}_0}}{\int_S K(z)e^{-2\hat{w}_0} dA}\right) = 0
   \vspace{5pt} \\
   \Jcal(-2\hat{w}_0) = \inf\limits_E\Jcal;
     \end{array} \right.
  \qquad
\ene
%Note that since $g=2$, the hyperbolic area $|S| = 4\pi$.
\item (Concentration) \underline{or}, the functional $\Jcal$ (with $K=\|\alpha\|_\sigma^2$) does not attain its infimum in $E$, and along a 
sequence $t_n\to 0$, there holds:
\be\label{3-1}
t_n^2K(z)e^{-2\tilde{u}_{t_n}} \rightharpoonup 4\pi\delta_{p_0},
\ene
weakly in the sense of measure, for some $p_0\in S$ such that $\alpha(p_0) \not= 0$,
\be\label{3-2}
(\tilde{u}_{t_n} - \fint_S \tilde{u}_{t_n}) \to 4\pi G(\cdot, p_0), \ \text{ in }\ W^{1,q}(S), 1<q<2,
\ene
and uniformly in $C_{loc}^{2,\beta}(S\backslash\{p_0\}), 0<\beta<1$, with $G(\cdot,p)$ the unique Green's function of the Laplace operator 
$\Delta$ on the hyperbolic surface $(S, g_\sigma)$, as defined in \eqref{green} below.
\een
\end{thx}
\vskip 0.1in
Hence, from Theorem C we see that, in contrast to higher genera, when $g=2$ the (unstable) {\mi} of $S$ does not survive the passage 
to the limit, as $t \to 0^{+}$. Indeed for the induced metric $\tilde{g}_t^0 = e^{2\tilde{u}_t}g_\sigma$, there holds:
\beq
|S|_{\tilde{g}_t^0 } \to 0, \ \text{as } t \to 0^{+}.
\eeq 
Nonetheless, when the functional $\Jcal$ in \eqref{J} with $K = \|\alpha\|_\sigma^2$ attains its minimum in $E$, then we can still find, as 
$t \to 0^{+}$, a ``limiting" configuration for the ``blown-up" surface $(S, \hat{g}_t^0)$ with
\beq
\hat{g}_t^0 = e^{2\hat{u}_t}g_\sigma, \ \text{and }\ \ \hat{u}_t = \tilde{u}_t - \log t.
\eeq 
Indeed, along a sequence $t = t_n\to 0^{+}$, the blown-up sequence $(S, \hat{g}_{t_n}^0)$ converges (in the sense of Gromov-Hausdorff) to 
a surface $(S, \hat{g}^0)$ conformally equivalent to $(S,g_\sigma)$ with same total negative curvature $-4\pi$. It is interesting to note that 
$(S, \hat{g}^0)$ admits a non-positive curvature vanishing exactly as $|\alpha|$.

More generally, for any sequence of unstable solutions $u_n$ of \eqref{t-ge} with $t = t_n\to 0^{+}$, as $n\to +\infty$, we shall carry out an 
analogous blow-up analysis whose details are contained in Theorem D of Section \S 3. Roughly speaking, to any such sequence, we shall 
associate a divisor $D$ in $S$ of the following type:
\be\label{div}
D = 2\sum_{j=1}^m(1+n(p_j))\delta_{p_j}
\ene
with suitable $p_j \in S$ and $n(p_j) = 0$ if $\alpha(p_j)\not=0$, while for $p_j$ a zero of $\alpha$ (i.e. $\alpha(p_j) = 0$), then $n(p_j)$ is given 
by the corresponding multiplicity. Furthermore, 
\be\label{div2}
\chi(S) + |D| \le 0.
\ene
In case $\chi(S) + |D| < 0$, then as before, we find that the sequence of corresponding {\mi}s of $S$ with data $(\sigma,t_n\alpha)$ converges 
(in the sense of Gromov-Hausdorff) to a ``limiting" totally geodesic immersion into a three-dimensional hyperbolic cone-manifold which induces on 
$S$ a {\hym} with conical singularities along the divisor $D$ in \eqref{div}. While in case $\chi(S) + |D| =0$, then no such singular metric exist on 
$S$ and exactly as above, the given (unstable) {\mi}s do not survive the passage to the limit. Thus, as before, the asymptotic behavior of the 
blow-up sequence $u_n$ is described by a ``concentration-compactness" alternative, see Theorem D in section \S 3 for details.

It is an interesting open problem to see whether unstable solutions of \eqref{t-ge} can be constructed in such a way that their blow-up behavior 
matches a prescribed divisor of the type \eqref{div} and \eqref{div2}.

In conclusion, let us make a few remarks.

\br
It is well known that any $\alpha\in Q(\sigma)$ admits $4(g-1)$ zeroes, counting multiplicity. Seen from Theorems ~\ref{main2} and 
~\ref{main3}, we have that the blow-up of the unstable solution $\tilde{u}_t$, as $t\to 0$, cannot occur around a zero of $\alpha$. A more 
precise characterization of the blow-up point $p_0$ in Theorems ~\ref{main2} and ~\ref{main3} will be given in the sections \S4 and \S5. 
\er
\br
It is interesting to record that for $g=2$, the behavior of the unstable solution $\tilde{u}_t$, as $t\to 0^{+}$ depends on whether the 
Moser-Trudinger functional $\Jcal$ in \eqref{J} with weight $K=\|\alpha\|_\sigma^2$ attains in its infimum in $E$. Actually, exactly when 
$K(z) = \|\alpha\|_\sigma^2$, the existence of extrema for $\Jcal$ appears to be a delicate open problem. Indeed we shall see in section \S 5.2 
that in such case the functional $\Jcal$ just misses to satisfy the condition provided (for a general weight function $K$) in Theorem 7.2 of 
\cite{DJLW97}] which is sufficient to ensure the existence of a global minimum for $\Jcal$. 
\er
\br
From the blow-up analysis in Theorem D where Theorems ~\ref{main2} and ~\ref{main3} enter as special cases, we shall obtain a compactness 
result for solutions of \eqref{t-ge}. This will enable us to obtain a {\mi} of $S$ with prescribed total extrinsic curvature, 
$\rho = \int_Sdet_{g_0}(\text{II}) dA(g_0) \in (0, 4\pi (g-1))$ and data $(\sigma, t_\rho\alpha)$, with suitable $t_\rho \in (0,\tau_0)$. It would be 
interesting to investigate the dependence of $t_\rho$ {\wrt} $\rho$.
\er
\br
(Added to the proof) The statement of Theorem A and parts of its proof are somewhat similar to a result of Ding-Liu (\cite{DL95}) where they analyzed 
a problem of prescribing Gaussian curvature on closed surfaces. In fact from the variational point of view, both problems (with a parameter) admit 
similar structures. More precisely, we have the presence of a (strict) local minimum and a ``mountain pass geometry", for a suitable sharp range of the 
parameter involved. Typically, in this situation one can claim the existence of a stable and unstable solution (as in Theorem A and in \cite{DL95}) with 
the stable solution ``bifurcating" out of a known (trivial) solution of the problem at a limiting value of the parameter. Much more interesting is the description 
of the asymptotic behavior of the unstable solution, which reflects the particular nature of the geometrical problem in hands. This is the purpose of our 
Theorems B, C and D. For the Gauss curvature problem treated in \cite{DL95}, this goal has been pursued in \cite{BGS15}, and more recently by Struwe in 
\cite{Str20}, who obtained complete results. In fact, Struwe's result encourages the possibility to also complete our blow-up analysis in Theorem D when
blow-up occurs at a zero of the quadratic differential $\alpha$, see Remark 6.3 and 6.4 below. We thank the referee for pointing out these references. 
\er
%------------------------------------------------------
\subsection*{Plan of the rest of the paper:} In \S 2, we will provide several estimates before we move to prove Theorem ~\ref{main1} in section 
\S 4. Detailed blow-up analysis is conducted in sections \S 3 where we obtain Theorem D, and in \S5, where we prove Theorems ~\ref{main2} and 
~\ref{main3}. In \S 6, we extend the program to explore this {\mi} problem when prescribing the total extrinsic curvature.

%------------------------------------------------------
\subsection*{Acknowledgements} We wish to thank Biao Wang for his help with the figures. The research of Z.H. was supported by a grant from the 
Simons Foundation (\#359635, Zheng Huang). The research of M.L. was supported by MINECO grant MTM2017-84214-C2-1-P. The research of 
G.T. was supported by PRIN 2015 Project ``Variational methods with application to problems in mathematical physics and geometry" and MIUR 
excellence project, Department of Mathematics, University of Rome Tor Vargata CUP E83C18000100006.
%------------------------------------------------------

\section{Elementary estimates}
%------------------------------------------------------
%\subsection{Elementary estimates}
Before we proceed, we introduce more convenient notations. We set 
\beq
v = -2u
\eeq
and
\beq
K(z) = \|\alpha\|_\sigma^2=\frac{|\alpha|^2}{det(g_\sigma)},
\eeq
and rewrite the Gauss equation ~\eqref{t-ge} as follows:
\be\label{v-ge}
-\Delta v= 2t^2Ke^v -2(1 -e^{-v}),
\ene
where $v\in H^1(S)$, $t\ge 0$, and $K(z) \ge 0$ has finitely many zeroes, given by the zeroes of the prescribed {\hqd} $\alpha \in Q(\sigma)$, 
whose total number is $4g-4$, counting multiplicity.  
\bd
We call a function $v_t \in H^1(S)$ a solution of problem $(1)_t$ for $t\ge0$, if it solves the equation \eqref{v-ge}.
\ed
We collect some basic properties for the solutions of this problem.
\bl\label{elelem}
If $v$ is a solution of problem $(1)_t$, then we have 
\ben[{\bf i)}]
\item 
\be\label{gb}
t^2\int_S K(z)e^v dA + \int_Se^{-v}dA = 4\pi(g-1),
\ene
and in particular, 
\be\label{2-2-a}
(2\pi(g-1))^2 \ge t^2\int_S K(z)e^v dA\int_Se^{-v}dA,
\ene
\item $v\ge 0$ and $v\equiv 0$ if and only if $t=0$. Therefore $v(z) > 0$ for all $z\in S$ for any $t>0$.
\item If we write $v = w+ \ttc$, with $\int_Sw(z)dA =0$, and $\ttc =\fint_S v(z)dA$, then $\ttc > 0$ and 
\be\label{hat-v}
e^{\ttc} = \frac{2\pi(g-1)\pm\sqrt{(2\pi(g-1))^2-t^2\int_S K(z)e^w dA\int_Se^{-w}dA}}{t^2\int_S K(z)e^w dA}.
\ene
\een
\el
\bp
These properties follow by direct and elementary calculations. More specifically, to obtain \eqref{gb}, we integrate \eqref{v-ge} and apply the 
Gauss-Bonnet formula. At this point \eqref{2-2-a} is a direct consequence of \eqref{gb} and Schwarz inequality. 
\vskip 0.1in
In order to show (ii), we simply write $v = v^{+} - v^{-}$ where $v^{+}$ and $v^{-}$ are both non-negative. Using $v^{-}$ as a test function we have 
\beq
\int_S \nabla v\nabla v^{-}\ dA = 2\left(\int_S t^2Ke^vv^{-} + (e^{-v}-1)v^{-}\ dA\right),
\eeq
and this is equivalent to 
\beq
-\int_S |\nabla v^{-}|^2\ dA = 2\left(\int_S t^2Ke^vv^{-} + (e^{v^{-}}-1)v^{-}\ dA\right),
\eeq
Since $v^{-} \ge 0$, we find 
\beq
\int_S (e^{v^{-}}-1)v^{-}\ dA \le 0.
\eeq
But the scalar function $(e^x-1)x \ge 0$ whenever $x \ge 0$, we deduce that $v^{-} \equiv 0$ and therefore $v \ge 0$. Furthermore, since $v \in H^1(S)$, 
the righthand side of \eqref{2-2-a} is in $L^p$ for any $p \ge 1$. Standard regularity results for elliptic equations imply that $v$ is in fact smooth. Then the 
conclusion follows from the strong {\maxp}. 
\vskip 0.1in
Finally (iii) follows by using $v = w+ \ttc$ in \eqref{gb} and by solving {\wrt} $\ttc$.
\ep
For later use and according to the choice of the sign in \eqref{hat-v}, we set:
\be\label{c_pm}
\ttc_{\pm}= \log\{\frac{2\pi(g-1)\pm\sqrt{(2\pi(g-1))^2-t^2\int_S K(z)e^w dA\int_Se^{-w}dA}}{t^2\int_S K(z)e^w dA}\}.
\ene
%---------------
\bl (\cite{HL12})
If $(1)_t$ admits a solution, then 
\be\label{upper-t}
t \le {\frac{1}{\fint_S2\sqrt{K(z)}dA}} =t^{\ast}.
\ene
\el
%----------------
\bp
Inequality \eqref{upper-t} easily follows from \eqref{gb} and the Cauchy-Schwarz inequality, as seen in \cite{HL12}.
\ep
In particular, this lemma implies that, for $t>t^{\ast}$, there does not exist any {\mi} of $S$ with data $(\sigma, t\alpha)$.
\bl\label{2-4}
Let $v = w+ \ttc$ be a solution to problem $(1)_t$, where $w$ and $\ttc$ are as in Lemma ~\ref{elelem}. We have 
\ben[{\bf i)}]
\item For any $s \in [1,2)$, there exists a constant $C_s > 0$ such that
\be\label{bd-W-Ls}
\|\nabla w\|_{L^s} \le C_s.
\ene
\item There exists a constant $C$ (independent of $t$) such that:
\be\label{2-8}
\int_S Ke^w\ dA \ge C.
\ene
\item For any small $\delta>0$, there exists a constant $C_\delta > 0$ (independent of $t$) such that,
\be\label{2-9}
0 < \ttc \le C_\delta,
\ene
for any $t \ge \delta$.
\een
\el
\bp
Since the righthand side of \eqref{v-ge} is {\ub} in $L^1$-norm because of \eqref{gb}, independent of $t$, then \eqref{bd-W-Ls} follows from Stampacchia 
elliptic estimates. In particular, for any $p\ge 1$, there exists a constant $C_p>0$ such that $\|w\|_p \le C_p$. Denote by 
$\kappa_0 = \frac{\|K\|_\infty}{\|K\|_1}$, then by Jensen's inequality we find:
\bear
\int_SK(z)e^wdA &\ge& \left(\int_SK\ dA\right)e^{\int_S\frac{K}{\|K\|_1}w\ dA} \\
&\ge& \|K\|_1e^{-\kappa_0\|w\|_1} \\
&\ge& e^{-C}, 
\eear
and \eqref{2-8} follows.

From \eqref{c_pm}, for $t \ge \delta$, we have:
\bear
0\le \ttc &\le& \log(\frac{4\pi(g-1)}{t^2\int_SK(z)e^wdA}) \\
&\le& \log(\frac{4\pi(g-1)}{\delta^2}) -\log(\int_SK(z)e^wdA),
\eear
and \eqref{2-9} follows from \eqref{2-8}.
\ep
In the next section we complete the information of Lemma ~\ref{2-4} by means of a more accurate blow-up analysis.
%------------------------------------------------------
%\subsection{Green's Function}

%%%%%%%%%%%%%%%%%%%%%%%%%
\section{Blow-up Analysis} 
In this section we will study in details the general asymptotic behavior of a blow-up sequence of solution for problem $(1)_t$. As an application, we 
will prove Theorems B and C in sections \S 5 and \S 6.

Let us denote by $G(q,p)$ the Green's function (for the hyperbolic Laplace operator) defined as follows:
\be\label{green} 
   \left\{
  \begin{array}{cl}
    -\Delta G = \delta_p - \frac{1}{4\pi(g-1)}   
    \vspace{5pt} \\
    \int_SG(q,p)dA(q)=0.                     
  \end{array} \right.
  \qquad
\ene
Here $\delta_p$ is the Dirac measure with pole at the point $p \in S$. It is well-known (see for instance \cite{Aub98}) that 
\beq
G(p,q) = G(q,p), \text{ for } p\not= q,
\eeq
and
\be\label{verde}
G(q,p) = -\frac{1}{2\pi}\log(dist(q,p)) + \gamma(q,p),
\ene 
where $\gamma \in C^{\infty}(S\times S)$ is the regular part of the Green's function $G$. 

%%%%%%%%%%%%%%%%%%%%%%%%%%%%%%%%%%%
\subsection{Mean field formulation}
Our blow-up analysis about solutions of problem $(1)_t$ is based on well-known results concerning blow-up solutions of Liouville type 
equations in mean field form (see \cite{BM91, LS94, BT02}). Thus, we proceed first to reformulate problem $(1)_t$ to a mean field type 
equation. To this end, we let $v$ be a solution of problem $(1)_{t}$ and set: 
\be\label{rho0}
\rho = t^2\int_S Ke^{v}dA.
\ene
By the conservation identity \eqref{gb}, we know that $\rho \in (0, 4\pi(g-1))$. As before (in Lemma ~\ref{elelem}), we set $v = w+ \ttc$, where 
$\int_Sw(z)dA =0$, and $\ttc =\fint_S v(z)dA$. 
\bl\label{mfe}
If $v= w+ \ttc$ is a solution of the problem $(1)_{t}$ satisfying \eqref{rho0} with $\rho \in (0, 4\pi(g-1))$, then $w$ satisfies:
\be\label{w} \left\{
  \begin{array}{cl}
    -\Delta w = 2\rho\left(\frac{K(z)e^{w}}{\int_S K(z)e^{w} dA}-\frac{1}{|S|} \right) +2(4\pi(g-1)-\rho)\left(\frac{e^{-w}}{\int_Se^{-w} dA}-\frac{1}{|S|} \right)
    \vspace{5pt} \\
    \int_Sw(z)dA =0,                    
  \end{array} \right.
  \qquad
\ene
Vice versa, if $w_\rho$ solves equation \eqref{w}, with $\rho \in (0, 4\pi(g-1))$, then by setting 
\beq
\ttc_\rho =\log\left(\frac{\int_Se^{-w} dA}{4\pi(g-1)-\rho} \right),
\eeq
and 
\beq
t_\rho^2 =\frac{\rho(4\pi(g-1)-\rho)}{(\int_S K(z)e^{w} dA)(\int_Se^{-w} dA)},
\eeq
we see that $v = w_\rho+\ttc_\rho$ is a solution to problem $(1)_{t_\rho}$.
\el
\bp
This follows by direct and simple calculations.
\ep
%%%%%%%%%%%%%%%%%%%%%%%%%%%%%%%%%%%%%
\subsection{General blow-up}  
Recall that $K(z) = \frac{|\alpha(z)|^2}{det(g_\sigma)}$. Denote by $\{q_1, \cdots, q_N\}$ the (finite) set of distinct zeroes of $\alpha$, i.e.,
\beq
\alpha(z) = 0 \iff z =q_i \text{ for some }i\in \{1,\cdots, N\},
\eeq 
and let $n_i$ be the multipilcity of $q_i$. It is well known that, $\sum\limits_{i=1}^Nn_i =4(g-1)$. 

Our main result in this section is the following theorem, which in particular indicates that blow-up occurs only when $t\to 0$.
\begin{thx}\label{main4}
Let $v_n$ be a solution of the problem $(1)_{t_n}$ such that 
\beq
\max_{S} v_n \to +\infty, \text{ as } n\to\infty,
\eeq
then, as $n\to\infty$, 
\be\label{401a}
t_n \to 0,
\ene
and
\be\label{m1}
t_n^2\int_S K(z)e^{v_n}\ dA \to 4\pi m, \ \ \ \text{ for some } m\in\{1, \cdots, g-1\},
\ene
where $g\ge 2$ is the genus of $S$. Furthermore,
\ben [{\bf (i)}]
\item
if $1 \le m < g-1$, then ({\aas}),
\ben
\item
there exist $\{p_1, \cdots, p_s \}\subset S$ (called \underline{blow-up points}), and sequences $\{p_{j,n}\}\subset S$ such that $p_{j,n}\to p_j$ 
and $v_n(p_{j,n}) \to +\infty$ as $n\to+\infty, \ j =1, \cdots, s$. Moreover,
\be\label{1-c}
t_n^2K(z)e^{v_n} \rightharpoonup 4\pi\sum_{j=1}^s( 1+n(p_j))\delta_{p_j}
\ene
weakly in the sense of measure, with
\be\label{1-d}
n(p_j) =   \left\{
  \begin{array}{cl}
    &0, \ \ \ \ \ \text{ if }\alpha(p_j)\not= 0
    \vspace{5pt} \\
    &n_i,  \ \ \ \ \text{ if }\alpha(p_j)= 0 \text{ for some }1\le j\le s, \text{ and } p_j = q_i,                    
  \end{array} \right.
  \qquad
\ene
and \ $m = \sum\limits_{j=1}^s( 1+n(p_j))$.
\item
$v_n \rightharpoonup v_0$ weakly in $W^{1,q}(S)$, for some $v_0$ and $1< q<2$, and uniformly in 
$C_{loc}^{2,\beta}(S\backslash\{p_1, \cdots, p_s\}$ for some $\beta \in (0,1)$. We also have 
\be\label{4.8a}
e^{-v_n}\to e^{-v_0},
\ene
strongly in $L^p(S)$ for all $p \ge 1$, and $v_0$ is the unique solution of the following equation on $S$:
\be\label{1-e}
-\Delta v_0 = 2\left(4\pi\sum_{j=1}^s( 1+n(p_j))\delta_{p_j} + e^{-v_0} -1\right).
\ene
\een
\vskip 0.1in
\item
If $m = g-1$, then by setting $v_n = w_n+ \ttc_n$, where $\int_Sw_n(z)dA =0$, and $\ttc_n =\fint_S v_n(z)dA$, we have ({\aas}):
\be\label{1-f}
\ttc_n \to +\infty, \text{ and } w_n \rightharpoonup w_0 \text{ weakly in } W^{1,q}(S), 1<q<2.
\ene
Furthermore we have the following alternatives:
\ben[{\bf (a)}]
\item (Compactness) \underline{either},
\beq
w_n \to w_0, \ \ \ \text{strongly in }\ H^1(S), \ \text{and in any other relevant norm},
\eeq
and 
\beq
t_n^2K(z)e^{v_n} \to \frac{4\pi(g-1)K(z)e^{w_0}}{\int_S K(z)e^{w_0}\ dA} \ \ \ \text{ strongly in }H^1(S),
\eeq
with $w_0$ satisfying the following equation on $S$:
\beq\left\{
  \begin{array}{cl}
    -\Delta w_0 = 8\pi(g-1) \left(\frac{K(z)e^{w_0}}{\int_S K(z)e^{w_0} dA} -\frac{1}{|S|}\right) 
    \vspace{5pt} \\
    \int_Sw_0(z)dA =0,                    
  \end{array} \right.
  \qquad
\eeq
\item 
 (Concentration) \underline{or}, for suitable points $\{p_1, \cdots, p_s\} \subset S$ (blow-up points), we have sequences 
$\{p_{i,n}\}\subset S$: $p_{i,n}\to p_i$ such that:
\beq
w_n(p_{i,n})=v_n(p_{i,n})-\fint_Sv_n \to +\infty, \ \ \text{as } n \to +\infty,
\eeq
and:
\beq
t_n^2K(z)e^{v_n} \rightharpoonup 4\pi\sum_{i=1}^s( 1+n(p_i))\delta_{p_i},
\eeq
weakly in the sense of measure, where $n(p_i)$ is defined in \eqref{1-d} with $\sum_{i=1}^s( 1+n(p_i)) = g-1$, and $w_0(z)$ satisfying:
\be\label{1-g}
w_0(z) = 8\pi\sum_{i=1}^s( 1+n(p_i))G(z, p), 
\ene
with $G(z,p)$ the unique Green's function in \eqref{green}. The convergence is uniform in $C_{loc}^{2,\beta}(S\backslash\{p_1, \cdots, p_s\})$.
\een
\een
\end{thx}
%---------------------
\bp 
As in the statement, we write $v_n = w_n+ \ttc_n$, where $\int_Sw_n(z)dA =0$, and $\ttc_n =\fint_S v_n(z)dA$. By Lemma 
~\ref{2-4}, we can always assume that, \aas, 
\beq
w_n \rightharpoonup w_0 \text{ weakly in } W^{1,q}(S),\  1<q<2. 
\eeq
Let 
\be\label{rho}
\rho_n = t_n^2\int_S Ke^{v_n}dA \in (0, 4\pi(g-1)).
\ene
So we can write $\int_Se^{-v_n}dA = 4\pi(g-1) - \rho_n$ (recall \eqref{gb}).

Denote by $\zeta_n$ the unique solution to the problem:
\be\label{zeta} \left\{
  \begin{array}{cl}
    -\Delta \zeta_n = 2(4\pi(g-1)-\rho_n)\left(\frac{e^{-w_n}}{\int_Se^{-w_n}}-\frac{1}{|S|} \right)\ \ \ \text{ on }S,
    \vspace{5pt} \\
    \int_S\zeta_n(z)dA =0.                  
  \end{array} \right.
  \qquad
\ene
Recall that $v_n>0$ in $S$, and 
\beq
(4\pi(g-1)-\rho_n)\frac{e^{-w_n}}{\int_Se^{-w_n}} = e^{-v_n}.
\eeq
Therefore, we know that the righthand side of \eqref{zeta} is {\ub} in $L^\infty(S)$. Thus, by elliptic estimates, we derive that $\zeta_n$ is 
{\ub} in $C^{2,\beta}(S)$-norm, with $\beta \in (0,1)$. So, {\aas}, we can assume that,
\be\label{zeta2}
\zeta_n \to \zeta_0 \ \ \ \text{ in }C^2(S)\text{-norm}, \ \ \text{as } n\to +\infty.
\ene
We define, somewhat abusing our notation to use $z_n$ as follows, 
\be\label{z1}
z_n = w_n - \zeta_n,
\ene
which satisfies the following mean field type equation:
\be\label{6*} \left\{
  \begin{array}{cl}
    -\Delta z_n = 2\rho_n\left(\frac{K_n(z)e^{z_n}}{\int_SK_n(z)e^{z_n}}-\frac{1}{|S|}\right)\ \ \ \text{ on }S,
    \vspace{5pt} \\
    \int_Sz_n(z)dA =0,               
  \end{array} \right.
  \qquad
\ene
with 
\be\label{4.16a}
K_n = K e^{\zeta_n} \to Ke^{\zeta_0} \ \ \ \text{in }C^2(S), \ \ \text{as } n\to +\infty.
\ene
At this point, we recall the following well-known ``concentration-compactness" result of \cite{BM91, LS94, BT02} which we state in a form suitable 
for our situation: 
%--------------------
\bt\label{star} (\cite{BM91, LS94, BT02})
Assume \eqref{6*} with $\rho_n\in (0, 4\pi(g-1))$ and \eqref{4.16a} with $K(z) = \frac{|\alpha|^2}{det(g_\sigma)}$, then, {\aas}, as $n\to +\infty$:
\beq
z_n \rightharpoonup z_0 \text{ weakly in } W^{1,q}(S), 1<q<2,
\eeq
and the following alternative holds:
\ben[{\bf (1)}]
\item \underline{either}, $\max_S z_n \le C$, and {\aas}, as $n\to+\infty$:
\beq
\rho_n\to \rho_0\in [0, 4\pi(g-1)],
\eeq
\beq
z_n \to z_0 \ \ \ \text{ strongly in }H^1(S),
\eeq
and in any other relevant norm, with $z_0$ satisfying:
\be\label{z0}\left\{
  \begin{array}{cl}
    -\Delta z_0 = 2\rho_0\left(\frac{he^{z_0}}{\int_She^{z_0}}-\frac{1}{|S|}\right)\ \ \ \text{ on }S,
    \vspace{5pt} \\
    \int_Sz_0dA =0,               
  \end{array} \right.
  \qquad
\ene
with $h=Ke^{\zeta_0}$ (see \eqref{zeta2}).
\item 
\underline{or}, there exist (blow-up) points $\{p_1, \cdots, p_s\} \subset S$, and sequences $\{p_{j,n}\}\subset S$: $p_{j,n}\to p_j$, 
such that $z_n(p_{j,n})\to+\infty$ as $n \to+\infty$, such that,
\be\label{4.17c}
\rho_n \frac{K_n(z)e^{z_n}}{\int_SK_n(z)e^{z_n}} \rightharpoonup 4\pi\sum_{j=1}^s(1+n(p_j))\delta_{p_j},
\ene
weakly in the sense of measure, with $n(p_j)$ defined in \eqref{1-d}. In particular,
\be\label{7a}
\rho_n \to \rho_0=4\pi\sum_{j=1}^s(1+n(p_j)) \in 4\pi\{1,\cdots, (g-1) \}\subset 4\pi\mathbb{N},
\ene
\beq
z_n \to z_0 \ \ \text{uniformly in } C_{loc}^{2,\beta}(S\backslash\{p_1, \cdots, p_s\}),\ \ \  0<\beta<1,
\eeq
and 
\beq
z_0(x) = 8\pi\sum_{i=1}^s(1+n(p_i))G(x,p_i).
\eeq
\een
\et
\bp
See Theorem 5.7.65 in \cite{Tar08}.
\ep
%-----------
Back to the proof of Theorem D. We apply Theorem ~\ref{star} to $z_n$ in \eqref{z1}, and {\aas}, we have: 
\beq
z_n \rightharpoonup z_0 \text{ weakly in } W^{1,q}(S),\ 1<q<2.
\eeq
Recall that by assumption: $\max\limits_S v_n \to \infty$, and so in case alternative (1) in Theorem ~\ref{star} holds, in view of \eqref{zeta2} and 
\eqref{z1}, we see that necessarily
\beq
\ttc_n \to +\infty, \ \text{and }w_n \to w_0 \ \text{strongly }, \ \text{ as } n\to\infty.
\eeq
Therefore, by dominated convergence, we derive:
\be\label{rho4}
4\pi(g-1)-\rho_n = \int_Se^{-v_n}dA \to 0, \text{ as } n\to\infty.
\ene
So $\rho_n\to 4\pi(g-1)$, and by \eqref{rho} we deduce that \eqref{m1} must hold with $m = g-1$ in this case.

This covers the compactness part in the statement (ii). Furthermore, by part (iii) of Lemma ~\ref{2-4}, we see also that 
$t_n\to 0$, as $n\to\infty$, and so \eqref{401a} holds.

Next we assume that alternative (2) in Theorem ~\ref{star} holds. Then by virtue of \eqref{rho}, \eqref{4.17c}, and \eqref{7a}, we check that \eqref{m1} holds with 
$m:=\sum\limits_{j=1}^s(1+n(p_j))$. We consider first the case where $1\le m < g-1$. As a consequence, 
\beq
\int_Se^{-v_n} \to 4\pi(g-1-m) > 0,
\eeq 
and so necessarily $\ttc_n =\fint_S v_n$ must be {\ub}. Thus {\aas}, we have: $v_n \rightharpoonup v_0$ weakly in $W^{1,q}(S)$, 
$1< q<2$, and also uniformly in $C^{2,\beta}(S\backslash\{p_1, \cdots, p_s\})$ with $0<\beta<1$, for some $v_0$. Furthermore, by dominated 
convergence, we see that  
\beq
e^{-v_n} \to e^{-v_0} \ \text{ in }L^p(S),\ \forall\ p\ge 1.
\eeq
As a consequence $v_0$ satisfies \eqref{1-e}. In other words, we have verified that part (i) holds in this case.

Since,
\beq
v_n\to v_0, \ \text{and }\ t_n^2Ke^{v_n} \to 0,  \text{ as } n \to\infty,
\eeq
uniformly on compact sets of $S\backslash\{p_1, \cdots, p_s\}$, we may conclude that \eqref{401a} must hold as well.

Finally, when we have alternative $(2)$ with $m = g-1$, then necessarily 
\beq
\int_Se^{-v_n} \to 0,  \text{ as } n \to\infty.
\eeq 
As a consequence, $\ttc_n \to +\infty$, as $n\to\infty$ and this implies as above, that $t_n\to 0$, by part (c) of Lemma ~\ref{2-4}.
Thus we have verified \eqref{401a}, \eqref{1-f}, and \eqref{m1} with $m = g-1$. At this point, alternative (2) of Theorem~\ref{star} 
in this situation gives exactly the (concentration) part (b) of (ii).

Finally, since $w_0$ in \eqref{1-f} satisfies:
\be\label{w03}\left\{
  \begin{array}{cl}
    -\Delta w_0 = 8\pi\sum\limits_{j=1}^s(1+n(p_j))\delta_{p_j}-2=\sum\limits_{j=1}^s8\pi(1+n(p_j))(\delta_{p_j}-\frac{1}{|S|})
    \vspace{5pt} \\
    \int_Sw_0 =0,                    
  \end{array} \right.
  \qquad
\ene
we see that \eqref{1-g} holds, and the proof is complete.
 \ep
 %---------------
 Concerning the location of the (possible) blow-up points $\{p_1, \cdots, p_s\}$ of $v_n$, we can use well-known results (\cite{OS05}) 
 which apply to the blow-up points of the sequence $z_n$ in \eqref{6*}. Thus, according to Theorem 2.2 of (\cite{OS05}), we conclude 
 that, if $p_j$ is a blow-up point with $\alpha(p_j)\not=0$, then in conformal coordinates around $p_j$ there holds:
 \be\label{OS}
 \nabla_z\left(8\pi\gamma(z,p_i) +8\pi\sum_{j\not=i}(1+n(p_i))G(z,p_j)+\log h\right)|_{z=p_i} =0,
 \ene
 with $i\in\{1,\cdots,s\}$, and
 \be\label{h}
 h=Ke^{\zeta_0}, 
 \ene
 with $\zeta_0$ in \eqref{zeta2} and $K=\|\alpha\|_\sigma^2$.

 Notice in particular that when $m = g-1$, then the function $\zeta_0 \equiv 0$, and \eqref{OS} provides a well-known necessary condition for 
 blow-up at  $\{p_1, \cdots, p_s\}$. Indeed, in case of non-degeneracy, \eqref{OS} turns out to be also a sufficient condition for the construction 
 of blow-up solutions at $\{p_1, \cdots, p_s\}$ for mean field equations on surfaces, see for instance \cite{CL03}.

 On the contrary, when $1\le m<g-1$, the condition \eqref{OS} is more involved since the function $\zeta_0$ is nonzero and satisfies the following equation:
 \be\label{4.22} \left\{
  \begin{array}{cl}
    -\Delta \zeta_0 = 8\pi(g-1-m)\left(\frac{e^{-8\pi\sum\limits_{j=1}^s(1+n(p_j))G(z,p_j)}e^{-\zeta_0}}{\int_Se^{-8\pi\sum\limits_{j=1}^s(1+n(p_j))G(z,p_j)}e^{-\zeta_0}\ dA}-\frac{1}{|S|}\right)\ \ \ \text{ on }S,
    \vspace{5pt} \\
    \int_S\zeta_0(z)\ dA =0.               
  \end{array} \right.
  \qquad
\ene

 So $\zeta_0$ itself depends on the blow-up points $\{p_1, \cdots, p_s\}$. Therefore it would be interesting to see whether one can find a 
 (nondegenerate) set of points satisfying \eqref{OS}, \eqref{h} and \eqref{4.22} which turns out to be the blow-up set of a sequence of bubbling 
 solutions for $(1)_t$, along a sequence of $t$'s going to zero.

As a consequence of Theorem D, we know that blow-up can only occur as $t\to 0$. Therefore, we can complete the (uniform) estimates given in 
Lemma ~\ref{2-4} as follows:
\bcor\label{4-2}
For any $\delta > 0$, there exists a constant $C_\delta > 0$ such that any solution $v$ of the problem $(1)_t$ with $t\ge \delta$ satisfies:
\beq
\|v\|_\infty \le C_\delta.
\eeq 
\ecor

Actually, by means of elliptic estimates we know that the $L^\infty(S)-$norm above can be replaced by any other stronger norm. Theorem D can 
be better interpreted in terms of the mean field formulation of problem $(1)_t$, and gives 
the following ``compactness" result:
\bcor\label{4.4}
Let $w_n$ be a sequence of solutions for \eqref{w} with $\rho = \rho_n$, and $\rho_n \to \rho_0 \in (0, 4\pi(g-1)]\backslash\{4\pi m, 1\le m \le g-1\}$. 
Then {\aas}, $w_n\to w_0$ in $H^1(S)$ (and any other relevant norm), with $w_0$ a solution of \eqref{w} with $\rho = \rho_0$.
\ecor

\bcor\label{3.5}
For every compact set $\mathcal{A}\subset (0, 4\pi(g-1)]\backslash\{4\pi m, 1\le m \le g-1\}$, the set of solutions of \eqref{w} with 
$\rho\in\mathcal{A}$ is {\ub} in $C^{2,\beta}(S), \ 0<\beta<1$.
\ecor
%------------------------------------------------------
\section{Proof of theorem ~\ref{main1}}
%------------------------------------------------------

%------------------------------------------------------
In this section, we will prove parts of Theorem ~\ref{main1}, in various steps. In this way, we obtain a detailed description of the lower branch 
of the bifurcation solution curve $\Ccal$. We shall take advantage of the variational formulation of problem $(1)_t$. Indeed, it is easy to 
verify that (weak) solutions of problem $(1)_t$ correspond to the critical points of the following functional:
\be\label{I}
\Ical_t(v) = \frac12\|\nabla v\|_2^2 - 2t^2\int_SK(z)e^vdA + 2\int_Se^{-v}dA + 2\int_SvdA, \ \forall\ v\in H^1(S).
\ene

%%%%%%%%%%%%%%%%%%%%%%%%%%%%%%%%%%%%
\subsection{First bending point}
We define the following two sets: 
\be\label{4.2.a}
\Lambda = \{t\ge 0: (1)_t \text{ admits a solution } \},
\ene
and 
\be\label{4.2.b}
\Lambda_s = \{t\ge 0: (1)_t \text{ admits a {\it stable} solution } \}.
\ene
Clearly $\Lambda_s\subseteq\Lambda \subset [0,t^{\ast}]$, with $t^{\ast}$ given in \eqref{upper-t}. Furthermore, since the problem 
$(1)_{t=0}$ only admits the trivial solution $v =0$, which is strictly stable, we see that, $\Lambda_s$ is nonempty and, 
\be\label{ts} 
0 < \tau_0 = \sup\{\Lambda_s\} \le t_0= \sup\{\Lambda\}. 
\ene
We aim to show that $\Lambda = \Lambda_s$ and $t_0 = \tau_0$.

To this purpose, we observe firstly that, by the estimates in Corollary ~\ref{4-2} and a limiting argument, we know that problem $(1)_{\tau_0}$ 
admits a stable solution $v_0$ which is also degenerate. According to the language of Crandell and Rabinowitz (\cite{CR80}), $(v_0,\tau_0)$ 
defines a ``bending point" for the curve of solutions of problem $(1)_t$, given by the zero set of the map 
\beq
F(v,t) = \Delta v + 2- 2(e^{-v}+t^2Ke^v): C^{2,\beta}(S)\times \R\to C^{0,\beta}(S),
\eeq
with $0<\beta<1$.

To establish Theorem ~\ref{Uh}, Uhlenbeck (\cite{Uhl83}) showed that actually $\tau_0$ is the only value for which the problem $(1)_t$ 
admits a degenerate stable solution.

\bpo\label{bending}
The problem $(1)_t$ admits a degenerate stable solution only at $t = \tau_0$. Moreover, for any $t\in [0,\tau_0]$ problem $(1)_t$ admits a 
unique stable solution which forms a smooth monotone increasing curve ({\wrt} $t$), and it is strictly stable for $t\in [0,\tau_0)$.
\epo
\bp
Let $v_0$ be the degenerate stable solution for $(1)_{\tau_0}$. We know that $(v_{0},\tau_0)$ must correspond to a ``bending point" for the set: 
$F(v,t) =0$, around  $(v_{0},\tau_0)$ (see \cite{CR80}).In other words, for $\epsilon > 0$ small, there exists a smooth curve 
$(v(s), t(s))$ satisfying $F(v(s),t(s)) = 0$ for all $s \in (-\epsilon,\epsilon)$, such that: $v(0) = v_0, \ t(0) = \tau_0$, $\dot{t}(0) =0$, and 
$\dot{v} > 0$ (i.e. $v(s)$ is increasing), where we have used the ``dot" to denote derivatives {\wrt} $s$.

Uhlenbeck showed further that $\ddot{t} (0)<0$ (see \cite{Uhl83}), so that: 
\beq
t(s) < t_0, \ \forall\ s\in (-\epsilon,\epsilon)\backslash\{0\},
\eeq
and in particular, $\dot{t}(s) > 0$ for $s\in (-\epsilon,0)$. This implies that $v(s)$ is strictly stable for every $s\in (-\epsilon,0)$.

Note that the same local description would hold for any other (possible) degenerate stable solution for which the corresponding of \eqref{4.2.b} 
would hold. This shows that if we continue the lower branch $(v(s),t(s)), \ s\in (-\epsilon,0)$ with the Implicit Function Theorem, we see that it 
cannot join with another degenerate stable solution at lower value of $t$. Instead, the lower branch can be continued until it joins the trivial 
solution at $t=0$. Thus, we can conclude that there exists a smooth, increasing curve of strictly stable solutions of problem $(1)_t, \ t\in [0,\tau_0)$, 
which joins the trivial solution $v=0$ at $t=0$ with the degenerate stable solution $v_0$ at $t=\tau_0$. This also shows that for any 
$t\in [0,\tau_0]$, problem $(1)_t$ cannot admit any other stable solution. As otherwise we could argue as before to join  such a different solution to 
the trivial solution along another smooth curve of solutions, a contradiction to the non-denegeracy of the trivial solution $v\equiv 0$ at $t=0$. This 
concludes the proof.
\ep
Proposition ~\ref{bending} shows in particular that $\Lambda_s = [0,\tau_0]$ and it also furnishes a proof to Theorem ~\ref{Uh}.

In the next subsection, we shall prove, with the help of the sub/super solution method in variational guise (see \cite{Str00}), that problem $(1)_t$ 
admits a stable solution for any $t \in [0,t_0]$ where $t_0 =\sup\{\Lambda\}$. Consequently $\Lambda = \Lambda_s= [0. \tau_0]$.
%%%%%%%%%%%%%%%%%%%%%%%%%%%%%%%%
\subsection{Stable solutions}
We now prove the following theorem on stable solutions:
\bt\label{lower}
Let $t_0 = \sup\Lambda$ ($\Lambda$ in \eqref{4.2.a}). For any $t \in [0,t_0]$, problem $(1)_t$ admits a {\bf stable solution} $v_{1,t}$ which is 
strictly increasing {\wrt} $t$, and coincides with the smallest of the solutions (and supersolutions) of Problem $(1)_t$. Furthermore, for $t=t_0$, 
\beq
v_{1,t_0}(z) =\sup_{0\le t <t_0}v_{1,t}(z) < +\infty
\eeq
is the {\underline{unique}} solution for problem $(1)_{t_0}$. In particular: $t_0=\tau_0$ and $\Lambda_s=\Lambda = [0,\tau_0]$.
\et
\bp
Since for $t=0$, $v\equiv 0$ is the desired stable solution, we let $t \in (0,t_0)$ be fixed. By \eqref{ts}, we can find some $t_1 \in \Lambda$ such 
that, $0 < t < t_1$. We denote by $v_1> 0$ a solution for problem $(1)_{t_1}$, and observe that it defines a strict super-solution for problem $(1)_{t}$. 
Indeed, we have 
\beq
\int_S\nabla v_1\nabla \phi dA- 2t^2\int_SK(z)e^{v_1}\phi dA+ 2\int_Se^{-v_1}\phi dA + 2\int_Sv_1\phi dA > 0,
\eeq
for any $\phi\in H^1(S)$ with $\phi\ge 0$ a.e. in $S$, and $\phi\not\equiv 0$.

While $v_0 \equiv 0$ is an obvious strict sub-solution for problem $(1)_{t}$. We set
\be\label{z}
\Zcal = \{v\in H^1(S): 0 \le v \le v_1 \text{ a.e. in } S\}.
\ene
It is routine to verify that $\Zcal$ is a non-empty, convex and closed subset of $H^1(S)$. In addition, the functional $\Ical_t$ is bounded from below 
and coercive on $\Zcal$. Consequently the functional $\Ical_t$ attains its minimum value at a point $v_{1,t}$ in $\Zcal$, i.e., 
\be\label{v-t}
\Ical_t(v_{1,t}) = \min\limits_{\Zcal} \Ical_t, \ \ \ \text{and } 0 \le v_{1,t} \le v_1 \text{ on } S.
\ene
By the strict sub/super solution property of $v_0 \equiv 0$ and $v_1$ respectively, we have 
\beq
\Ical_t(v_{1,t}) < \min\{\Ical_t(0), \Ical_t(v_1) \},
\eeq
and therefore $v_{1,t}\not\equiv v_1$ and $v_t\not\equiv 0$. Following \cite{Str00}, we show next that $v_{1,t}$ is a {\cp} for $\Ical_t$, therefore a solution to 
$(1)_{t}$. We first let $\phi\in C^{\infty}(S)$ with $\phi \ge 0$, and for $\epsilon>0$ sufficiently small, we define 
\be\label{v-eps}
v_\epsilon = v_{1,t}+\epsilon\phi - \phi^\epsilon+\phi_\epsilon,
\ene
where
\beq
\phi^\epsilon = \max\{0,v_{1,t}+\epsilon\phi - v_1\} \ge 0,
\eeq
and 
\beq
\phi_\epsilon = \max\{0,-(v_{1,t}+\epsilon\phi) \} \ge 0.
\eeq
Therefore we have $v_\epsilon \in \Zcal$. By virtue of \eqref{v-t} and \eqref{v-eps}, we have:
\bea\label{ineq-inn}
0 &\le& \langle\Ical'_t(v_{1,t}), v_\epsilon-v_{1,t}\rangle \notag\\
&=& \epsilon\langle\Ical'_t(v_{1,t}),\phi\rangle - \langle\Ical'_t(v_{1,t}),\phi^\epsilon\rangle+ \langle\Ical'_t(v_{1,t}),\phi_\epsilon\rangle.
\eea
We define a set
\beq
\Omega_\epsilon = \{p \in S: v_{1,t}(p)+\epsilon\phi (p) \ge v_1(p) > v_{1,t}(p)\}.
\eeq
We observe that, $|\Omega_\epsilon|$ the measure of $\Omega_\epsilon$ goes to zero as $\epsilon \to 0$. Moreover,
\bear
-\langle\Ical'_t(v_{1,t}),\phi^\epsilon\rangle &=& -\langle\Ical'_t(v_1),\phi^\epsilon\rangle - \langle\Ical'_t(v_{1,t})-\Ical'_t(v_1),\phi^\epsilon\rangle \\
&\le & -\int_S \nabla(v_{1,t}-v_1)\nabla\phi^\epsilon + 2t^2\int_SK(e^{v_{1,t}}-e^{v_1})\phi^\epsilon \\
&&\ \ \ \ \ \ -2\int_S(e^{-v_{1,t}}-e^{-v_1})\phi^\epsilon-2\int_S(v_{1,t}-v_1)\phi^\epsilon\\
&\le & -\int_{S} \nabla(v_{1,t}-v_1)\nabla \phi^\epsilon -2\int_S(v_{1,t}-v_1)\phi^\epsilon\\
&= & -\int_{\Omega_\epsilon} \nabla(v_{1,t}-v_1)\nabla (v_{1,t}-v_1+\epsilon\phi)  -2\int_{\Omega_\epsilon}(v_{1,t}-v_1)\phi^\epsilon \\
&\le & -\int_{\Omega_\epsilon} \nabla(v_{1,t}-v_1)\nabla (v_{1,t}-v_1+\epsilon\phi) \\
&&\ \ \ \ \ +2\int_{\Omega_\epsilon} (e^{-v_{1,t}}-e^{-v_1})\phi^\epsilon +2\int_{\Omega_\epsilon} (e^{-v_{1,t}}-e^{-v_1})\phi^\epsilon\\
&\le &-\epsilon \int_{\Omega_\epsilon} \nabla(v_{1,t}-v_1)\nabla\phi + 2\epsilon\int_{\Omega_\epsilon} (v_1-v_{1,t})\phi \\
&=& o(\epsilon) \ \ \text{ as } \epsilon \to 0.
\eear
Similar calculations show that 
 \beq
\langle\Ical'_t(v_{1,t}),\phi_\epsilon\rangle = o(\epsilon), \ \ \text{ as } \epsilon \to 0.
\eeq
Applying \eqref{ineq-inn}, we find:
\beq
\langle\Ical'_t(v_{1,t}),\phi\rangle \ge 0.
\eeq
We obtain the reverse inequality by replacing $\phi$ with $-\phi$. Since $C^\infty(S)$ is dense in $H^1(S)$, by a density argument we have now proved:
\beq
\langle\Ical'_t(v_{1,t}),\phi\rangle = 0, \forall \phi\in H^1(S).
\eeq
This implies that $v_{1,t}$ is a solution to problem $(1)_t$. Note that $v_{1,t} \not\equiv 0$ and $v_{1,t}\not\equiv v_1$, so by the {\maxp}, we have:
\be
0 < v_{1,t}(z)<v_1(z), \forall z \in S.
\ene
This ensures that $v_{1,t}$ is a local minimum for the functional $\Ical_t$ in $C^1(S)$-norm.

We are going to show that $v_{1,t}$ is actually a local minimum for $\Ical_t$ in $H^1(S)$-norm, and hence a stable solution for problem $(1)_t$. To 
this purpose, we argue by contradiction. Suppose there exists 
$v_n \in H^1(S)$ such that $\Ical_t(v_n) < \Ical_t(v_{1,t})$, and $v_n\to v_{1,t}$ in $H^1(S)$. Letting 
\beq
\delta_n^2 = \frac12\|\nabla(v_{1,t}-v_n)\|_2^2 + \|v_{1,t}-v_n\|_2^2\to 0, \ \ \text{as } n\to +\infty,
\eeq 
we may assume, without loss of generality, that 
\beq
\Ical_t(v_n) = \min\limits_{v\in H^1(S)}\{\Ical_t(v): \frac12\|\nabla(v-v_{1,t})\|_2^2 + \|v-v_{1,t}\|_2^2 = \delta_n^2 \}.
\eeq
This enables us to apply the Lagrange multiplier method, and for suitable $\lambda_n\in\R$, we find that $v_n$ satisfies 
\be\label{Lag}
-\Delta v_n = 2t^2K(z)e^{v_n}+2e^{-v_n}-2+\lambda_n(-\Delta(v_n-v_{1,t})+2(v_n-v_{1,t}))
\ene
We set
\be\label{eta_n}
\eta_n=\frac{v_n-v_{1,t}}{\sqrt{\|\nabla(v_{1,t}-v_n)\|_2^2 + 2\|v_{1,t}-v_n\|_2^2}}.
\ene
Then we see that $\|\eta_n\|_{H^1} \le 1$. So {\aas}, $\eta_n$ converges weakly in $H^1(S)$ to some function $\eta \in H^1(S)$, as $n\to +\infty$. 
Moreover, for any $p\ge 1$ we also have 
\be\label{eta_p}
\|\eta_n-\eta\|_p \to 0.
\ene
By recalling that $\Ical'_t(v_{1,t})= 0$, we compute:
\bea\label{I-differ}
0 &>&\Ical_t(v_n)-\Ical_t(v_{1,t}) \notag \\
&=& \Ical_t(v_n)-\Ical_t(v_{1,t}) + \Ical'_t(v_{1,t})(v_{1,t}-v_n) \notag\\
&=& 2t^2\int_SKe^{v_{1,t}}(1+v_n-v_{1,t}-e^{v_n-v_{1,t}})dA  \notag \\
&&\ \ \ +2\int_Se^{-v_{1,t}}(e^{v_{1,t}-v_n}-1-(v_{1,t}-v_n))dA+\frac12\|\nabla(v_n-v_{1,t}) \|_2^2.
\eea
Since $v_n \to v_{1,t}$ in $H^1(S)$, we know that, for any $p \ge 1$,  
\beq
e^{v_n-v_{1,t}} \to 1, \text{ in } L^p(S).
\eeq
Therefore by view of \eqref{I-differ} we may conclude that
\be\label{eta2}
\int_SK(z)e^{v_{1,t}}\frac{(e^{v_n-v_{1,t}}-1-(v_n-v_{1,t}))}{\|\nabla(v_{1,t}-v_n)\|_2^2 + 2\|v_{1,t}-v_n\|_2^2} \to \frac12\int_SK(z)e^{v_{1,t}}\eta^2,
\ene
and 
\be\label{eta3}
\int_Se^{-v_{1,t}}\frac{(e^{v_n-v_{1,t}}-1-(v_n-v_{1,t}))}{\|\nabla(v_{1,t}-v_n)\|_2^2 + 2\|v_{1,t}-v_n\|_2^2} \to \frac12\int_Se^{-v_{1,t}}\eta^2.
\ene

We claim that $\eta\not= 0$. To see this, by contradiction we assume $\eta = 0$, then 
\beq
\|v_n-v_{1,t}\|_2 = o(\|\nabla(v_n-v_{1,t}) \|_2),
\eeq
as $n \to\infty$. From \eqref{I-differ}, \eqref{eta2} and \eqref{eta3}, we find:
\bea\label{I-differ2}
0 &\ge&\lim_{n\to\infty}\frac{\Ical_t(v_n)-\Ical_t(v_{1,t})}{\|\nabla(v_{1,t}-v_n)\|_2^2 + 2\|v_{1,t}-v_n\|_2^2} \notag\\
&=& \frac12- t^2\int_SK(z)e^{v_{1,t}}\eta^2 +\int_Se^{-v_{1,t}}\eta^2 \notag\\
&=&\frac12,
\eea
which is impossible. Therefore $\eta\not=0$.

Using \eqref{v-ge} for $v_{1,t}$ and \eqref{Lag}, we see that:
\bea\label{lamb}
(1-\lambda_n)(-\Delta(v_n-v_{1,t}) + 2(v_n-v_{1,t})) &=& 2t^2K(z)(e^{v_n}-e^{v_{1,t}})\notag  \\ 
&&+2(e^{-v_n}-e^{-v_{1,t}}) + 2(v_n-v_{1,t}).
\eea
As above we find, as $n \to\infty$, 
\be\label{eta4}
\int_SK(z)(e^{v_n}-e^{v_{1,t}})\frac{(v_n-v_{1,t})}{\|\nabla(v_{1,t}-v_n)\|_2^2 + 2\|v_{1,t}-v_n\|_2^2} \to \int_SK(z)e^{v_{1,t}}\eta^2,
\ene
and 
\be\label{eta5}
\int_S(e^{-v_n}-e^{-v_{1,t}})\frac{(v_n-v_{1,t})}{\|\nabla(v_{1,t}-v_n)\|_2^2 + 2\|v_{1,t}-v_n\|_2^2} \to -\int_Se^{-v_{1,t}}\eta^2.
\ene

Now we have: 
\be\label{3.18a}
\lim_{n\to\infty}(1-\lambda_n) = 2t^2\int_SK(z)e^{v_{1,t}}\eta^2 + 2\int_S(1-e^{-v_{1,t}})\eta^2 >0.
\ene
So by \eqref{lamb} and \eqref{3.18a}, we can use elliptic regularity theory to conclude that $v_n \in C^1(S)$. Furthermore, the 
righthand side of \eqref{lamb} converges to zero in $L^p(S)$, for $p > 1$, as $n\to +\infty$. Consequently, by \eqref{3.18a}, we can 
use again elliptic estimates to show that $(v_n-v_{1,t})\to 0$ in $C^1$-norm. This is a contradiction to the fact that $v_{1,t}$ is a local 
minimizer for the functional $\Ical_t$ in $C^1$-norm. Therefore $v_{1,t} \in H^1(S)$ is a local minimum of $\Ical_t$ and hence a {\it stable} 
solution for the problem $(1)_t$ with $t \in (0,t_0)$.

So far we have shown that $\forall t \in [0,t_0)$, problem $(1)_{t}$ admits a stable solution $v_{1,t}$ (and infinitely many supersolutions). 
As a consequence,
\beq
t_0=\tau_0=\sup\Lambda_s,
\eeq
and by Proposition ~\ref{bending} we know also that problem $(1)_{\tau_0}$ admits a unique stable (degenerate) solution $v_0$.

For $t\in (0,\tau_0)$, to show that $v_{1,t}$ is the smallest among all solutions (and supersolutions) of problem $(1)_{t}$, we define for $z\in S$,  
\be\label{v_t}
v_t(z) = \inf\{v(z): v \text{\ a solution or supersolution of problem } (1)_t \} \ge 0.
\ene
Clearly, $v_t(z)$ defines a supersolution of problem $(1)_{t}$, in the sense that the following holds:
\be\label{3.?}
\int_S\nabla v\nabla\phi -2t^2\int_SK(z)e^v\phi -2\int_Se^{-v}\phi + 2\int_S\phi \ge 0,
\ene
for any $\phi\in H^1(S)$ with $\phi\ge 0$. Since $t>0$, we have $v_t\not\equiv 0$. Moreover, \eqref{3.?} can never hold with a strict sign, as 
otherwise we would be in position to apply the sub/super solution method as above, and obtain a solution of problem $(1)_{t}$ which is smaller 
than $v_t$, in contradiction with \eqref{v_t}. Hence $v_t$ is a solution of problem $(1)_{t}$ which, by definition, is the smallest solution of 
$(1)_{t}$ and strictly increasing {\wrt} $t \in [0,\tau_0)$.

We claim: 
\be\label{equal}
v_t = v_{1,t}.
\ene
To establish this claim, it suffices to show that $v_t$ is stable, so that \eqref{equal} follows by the uniqueness of stable solutions in 
Proposition ~\ref{bending}.

To this purpose, for $t\in (0,\tau_0)$, we use $v_t$ as a supersolution to problem $(1)_{s}$, for $0< s<t$. As above, for $s\in (0,t)$, we obtain a stable 
solution $\bar{v}_s$ of problem $(1)_{s}$ satisfying $0 < \bar{v}_s < v_t$ in $S$. By taking a sequence $s_n \nearrow t$, then by dominated 
convergence and elliptic estimates, we see that $\bar{v}_{s_n} \to \bar{v}$ in $H^1(S)$, with $\bar{v}$ a stable solution of problem $(1)_{t}$ and 
$0 < \bar{v} \le v_t$. Since $v_t$ is the smallest solution to $(1)_{t}$, we conclude that $\bar{v} \equiv v_t$, and so $v_t$ is stable and \eqref{equal} 
is established.

From \eqref{equal}, it also follows that, $\forall t\in [0,\tau_0)$, $v_{1,t}< v_0$, with $v_0$ the unique stable solution of problem $(1)_{\tau_0}$. 
So by the monotonicity property of $v_{1,t}$ in $t$, we find:
\beq
\lim_{t\nearrow \tau_0}v_{1,t}(z) = \sup_{0\le t<\tau_0}\{ v_{1,t}(z)\} = v_0(z), \text{ as }t\nearrow \tau_0,
\eeq
where again by dominated convergence and elliptic estimates, the convergence actually occurs uniformly in $C^{2,\beta}(S), 0<\beta<1$. Clearly, 
$v_0$ must define the smallest solution of problem $(1)_{\tau_0}$, i.e., $v_0 = v_{1,\tau_0}$. In fact we show that actually $v_0$ is the only solution 
of problem $(1)_{\tau_0}$.

To this purpose we argue by contradiction and assume there is another solution $v'$ for the problem $(1)_{\tau_0}$. By construction, $v_0$ is the 
smallest solution at $t=\tau_0$, so necessarily $v'>v_0$ on $S$. As seen in Proposition ~\ref{bending}, around $(v_0,\tau_0)$, we find a 
solution curve $(v(s), t(s))$ such that for $s \in (0,\epsilon)$ and $\epsilon >0$ sufficiently small, we obtain a solution $v(s)$ for the 
problem $(1)_{t(s)}$ such that $t(s) < \tau_0$ and $v_0<v(s)< v'$. Thus for $t\in (t(s), \tau_0)$ we find $v(s)$ as subsolution and $v'$ as 
supersolution for problem $(1)_{t}$.  So we can use the sub/super solution method again, and for $t \in (t(s),\tau_0)$, we obtain a stable solution 
for the problem $(1)_{t}$ which will be greater than $v_0$, and therefore greater than the smallest solution $v_{1,t}$. This is impossible, since the 
smallest solution is also the only (strictly) stable solution of problem $(1)_t$.
\ep
%------------------------------------------------------
\subsection{Compactness for the functional $\Ical_t$}
We will complete the proof of Theorem ~\ref{main1} in this subsection, namely we will prove the existence of an additional solution for each 
$t \in (0,t_0)$. This will extend a multiplicity result in (\cite{HL12}).

By combining Proposition ~\ref{bending} and Theorem ~\ref{lower}, we have now established that, $\forall \ t \in (0,\tau_0)$, $v_{1,t}$ is a strict local 
minimum for $\Ical_t(v)$ in $H^1(S)$. Furthermore, one checks that, for every $t \in (0,\tau_0)$:
\beq
\Ical_t(v_{1,t}+C) \to -\infty,
\eeq
as $C \to +\infty$. In other words, for $t \in (0,\tau_0)$,  the functional $\Ical_t$ admits a ``mountain-pass" structure 
(\cite{AR73}). Next we establish the following Palais-Smale (compactness) condition:
%%%%%
\bt\label{PS}
Suppose a sequence $\{v_n\} \in H^1(S)$ satisfies that, $\Ical_t(v_n) \to c$ and $\Ical'_t(v_n)\to 0$ as $n\to \infty$, then 
$\{v_n\}$ admits a convergent subsequence. In particular, $c$ is a critical value for the functional $\Ical_t$.
\et
%------------------------------
\bp
We show first that $v_n$ is uniformly bounded in $H^1$-norm. As in Lemma ~\ref{elelem}, we write $v_n = w_n+\ttc_n$, with 
$\int_Sw_n(z)dA =0$, and $\ttc_n =\fint_S v_n(z)dA$. Then, we have 
\be\label{ps1}
\langle \Ical'_t(v_n), 1\rangle = -2t^2e^{\ttc_n}\int_SK(z)e^{w_n}dA -2e^{-\ttc_n}\int_Se^{-w_n}dA+8\pi(g-1).
\ene
By assumption, $\langle \Ical'_t(v_n), 1\rangle = o(1)$ as $n\to\infty$. Applying Jensen's inequality, we have, as $n\to \infty$:
\beq
e^{-\ttc_n} \le e^{-\ttc_n}\fint_Se^{-w_n}dA \le 1 + o(1).
\eeq
Therefore for some suitable constant $C_0>0$, we find $\ttc_n \ge -C_0$. Now from \eqref{I} and \eqref{ps1}, we obtain, as $n\to \infty$: 
\be\label{ps2}
\Ical(v_n) = \frac12\|\nabla w_n\|_2^2 +4e^{-\ttc_n}\int_Se^{-w_n}dA -8\pi(g-1)+8\pi(g-1)\ttc_n+o(1).
\ene
By assumption, $\Ical(v_n)$ is uniformly bounded, so from \eqref{ps2} we also see that $\ttc_n$ is bounded from above, and that 
$\|\nabla w_n\|_2$ is {\ub}.

In conclusion we have $\|v_n\|_{H^1} \le C$ for some suitable $C>0$. Therefore along a subsequence, $v_n$ converges to some 
$v\in H^1(S)$ weakly. The convergence is strong in $L^p(S)$ for $p \ge 1$. In particular we have $\ttc_n \to \fint_SvdA$ and $\Ical'(v) = 0$. 
By the {\mt}, we also have:
\beq
\|e^{|v_n|}\|_{L^p} \le C_p, \forall p\ge 1.
\eeq 
In particular, for $p=2$, and as $n\to \infty$, 
\bear
o(1) &=& \langle \Ical'(v_n), v_n-v\rangle \\
&=& \langle \Ical'(v_n) - \Ical'(v) , v_n-v\rangle \\
&=& \|\nabla (v_n-v)\|_2^2 - 2t^2\int K(e^{v_n}-e^v)(v_n-v)-2\int_S(e^{-v_n}-e^{-v})(v_n-v)\\
&\ge&  \|\nabla (v_n-v)\|_2^2 - C\|v_n-v\|_2^2 \\
&=&  \|\nabla (v_n-v)\|_2^2 +o(1).
\eear
In other words, we have 
\beq
\|\nabla (v_n-v)\|_2 \to 0,
\eeq 
as $n\to \infty$. This completes the proof.
\ep
We can now apply the {\mtps} construction of Ambrosetti-Rabinowitz (\cite{AR73}) to obtain a second (unstable) mountain pass solution 
$v_{2,t} > v_{1,t}$ for all $t\in (0,t_0)$, 
satisfying:
\be\label{mtp1}
\Ical_t(v_{2,t}) =\inf_{\Gamma\in \Pcal}\max_{s\in [0,1]}\Ical_t(\Gamma(s)) > \Ical_t(v_{1,t}),
\ene
with the path space 
\be\label{pathspace}
\Pcal = \{\Gamma: [0,1]\to H^1(S) \text{ is continuous with } \Gamma(0)=v_{1,t},\  \Ical_t(\Gamma(1))\le \Ical_t(v_{1,t}) -10\}.
\ene
Clearly $\Pcal$ is not empty, since for $A>0$ sufficiently large we easily check that $\Gamma(s) = v_{1,t}+sA, \ s\in [0,1]$ lies in $\Pcal$.

Finally we show that the unstable solution $v_{2,t}$ will not stay bounded as $t\to 0$:
\bpo\label{unstable}
For $t \in (0,\tau_0)$, let $v_{2,t}$ be the {\mtps} solution obtained above. Then: 
\be\label{behavior}
\max\limits_{S} v_{2,t} \to +\infty, \text{ as } t\to 0.
\ene
\epo
\bp
We argue by contradiction. Suppose that, along a sequence $t_n\to 0$, we have 
\beq
0\le v_{2,t_n} \le C,
\eeq
for suitable constant $C>0$. Then by elliptic estimates ({\aas}), we find that $\{v_{2,t_n}\}$ converges strongly in $C^{2,\beta}(S)$ norm 
to $v\equiv 0$, the unique solution of problem $(1)_{t=0}$. But this is impossible, since for $t>0$ small, the stable solution $v_{1,t}$ is the 
only solution of $(1)_t$ contained in a small ball centered at the origin .
\ep
Clearly, by a similar argument, we see that any family of unstable solutions of \eqref{t-ge} admits the same blow-up behavior in \eqref{behavior}, 
as $t\to 0^{+}$.
%To proceed further, we provide next a general blow-up analysis for (sequence of) solutions for problem $(1)_t$.

%%%%%%%%%%%%%%%%%%%%%%%%%%%%%%%%%%%%%%%%%%
\section{Blow-up analysis: applications to mountain pass solutions}
In this section, we apply the general blow-up analysis of \S 3 to the {\mpsl} $v_{2,t}$ of problem $(1)_t$ obtained in Theorem A. The asymptotic 
behavior of $v_{2,t}$ differs when the surface has genus two or higher.

%%%%%%%%%%%%

By Proposition ~\ref{unstable} and Theorem D, we know that:
\be\label{5.1}
\liminf_{t\to 0^{+}}\left(t^2\int_SK(z)e^{v_{2,t}}dA\right) = 4\pi m, 
\ene 
for suitable $m\in\mathbb{N}$ satisfying $1\le m \le g-1$.

Our first goal is to prove that actually, $m=1$. We start with the case where the genus of the surface $S$ is at least three.
\subsection{Blow-up analysis when $g\ge 3$} 
%-----------------
\bt\label{ge3}
Let the genus $g \ge 3$. Then for $K = \|\alpha\|_\sigma^2$ we have:
\ben[{\bf (i)}]
\item
\be\label{sing1}
\lim_{t\to 0}t^2\int_SK(z)e^{v_{2,t}}dA = 4\pi.
\ene 
\item
As $t\to 0$,  
\beq
t^2Ke^{v_{2,t}}\to 4\pi\delta_{p_0},
\eeq
with some suitable $p_0 \in S$ such that $K(p_0) \not= 0$ (i.e. $\alpha(p_0)\not=0$).
\item
\beq
v_{2,t} \rightharpoonup v_0 \ \text{weakly in } W^{1,q}(S), \ 1<q<2,
\eeq
\beq
v_{2,t} \to v_0 \ \text{strongly in }  C_{loc}^{2,\beta}(S\backslash\{p_0\}), \ 0<\beta<1, 
 \eeq 
and, 
\beq
e^{-v_{2,t}}\to e^{-v_0} \ \text{strongly in } L^p(S), \ p\ge 1.
\eeq
Moreover, $v_0$ is the \underline{unique} solution to the following equation on $S$:
\beq
-\Delta v_0 = 8\pi\delta_{p_0}+2e^{-v_0} -2.
\eeq
\een
\et
%------------
In order to establish Theorem ~\ref{ge3}, we establish first the following estimates:
\bl\label{ap}
If the genus $g\ge 3$, then for a suitable constant $C>0$, we have:
\be\label{average}
0 \le \fint_Sv_{2,t}(z)dA \le C,
\ene 
and 
\be\label{ap2}
|\Ical_t(v_{2,t}) - 8\pi\log\frac{1}{t^2}|\le C, \ \ \forall t\in (0,\tau_0).
\ene 
\el
\bp
By virtue of Corollary ~\ref{4-2}, clearly it suffices to prove \eqref{average} and \eqref{ap2} as $t \searrow 0$. We start by showing:
\be\label{ap2a}
\Ical_t(v_{2,t})\le 8\pi\log\frac{1}{t^2} + C, \ \ \text{ as } t\searrow 0,
\ene
with a suitable constant $C>0$ (independent of $t$).

To this purpose we use sharp estimates obtained in \cite{DJLW97} in order to establish the existence of minimizers for the Moser-Trudinger 
functional $\Jcal$ in \eqref{J}. We fix $p \in S$ with $K(p) \not= 0$. As in \cite{DJLW97}, we use normal (polar) coordinates at $p$, 
centered at the origin, so that for $r =dist(q,p)$, we have:
\beq
8\pi G(r,\theta) = -4\log r+A(p)+b_1r\cos\theta+b_2r\sin\theta+\beta(r,\theta), \text{ as }r\to 0,
\eeq 
with $A(p) = 8\pi\gamma(p,p)$, suitable constants $b_1$ and $b_2$ depending on the {\hym} $g_\sigma$, and $\beta(r,\theta) = o(r)$ as $r\to 0$. 
Recall that $G(p,q)$ is the Green's function in \eqref{verde}, and $\gamma(p,q)$ its regular part.

We let $\eta$ be a standard cut-off function such that:
\beq\left\{
  \begin{array}{cl}
    \eta\in C_0^\infty(B_{2a_t}(p)),
    \vspace{5pt} \\
   \eta = 1 \text{ in } B_{a_t}(p),  
    \vspace{5pt}     \\
    \|\nabla\eta\|_{L^\infty}\le \frac{C}{a_t},           
  \end{array} \right.
  \qquad
  \eeq
  where $a_t>0$ is chosen in such a way that $a_t\to 0$ and $\alpha_t = \frac{a_t}{t}\to \infty$, as $t \to 0^{+}$.

  Now we let,
  \beq
  \varphi_t(r,\theta) =\left\{
  \begin{array}{cl}
    -2\log(r^2+t^2)+b_1\cos\theta+b_2\sin\theta, \ \ \ \text{ for } 0\le r \le a_t
    \vspace{5pt} \\
   8\pi G(r,\theta) - \eta\beta(r,\theta)-A(p)-2\log(1+\frac{1}{\alpha_t^2}), \ \ \ \text{ for } a_t< r \le 2a_t
    \vspace{5pt}     \\
     8\pi G(r,\theta) -A(p)-2\log(1+\frac{1}{\alpha_t^2}), \ \ \ \text{ for } 2a_t< r.         
  \end{array} \right.
  \qquad
  \eeq
 For $\varphi_t$, we can use well-known estimates. For example from the much sharper estimates derived in \cite{DJLW97} that we apply with 
 $\epsilon = t^2$, $\phi_\epsilon = \varphi_t+\log t^2$, and $\alpha =\alpha_t$, we obtain that, as $t\to 0$,
 \be\label{5-1}
 \int_S|\nabla\varphi_t|^2dA = 16\pi\log\frac{1}{t^2} - 16\pi+8\pi A(p) +o(1),
 \ene
  \be\label{5-2}
 \fint_S\varphi_t\ dA = - A(p) +o(1), \ \ \ \ \ \int_Se^{-\varphi_t}dA = O(1),
 \ene
 and 
   \be\label{5-3}
t^2\int_SK(z)e^{\varphi_t}dA = K(p)\pi+o(1).
 \ene
 Next we construct a suitable path in $\Pcal$ (defined in \eqref{pathspace}) as follows:
   \beq
  \Gamma_t(s) =\left\{
  \begin{array}{cl}
    (1-4s)v_{1.t}, \ \ \ \text{ for } 0\le s \le \frac14
    \vspace{5pt} \\
   (4s-1)\varphi_t, \ \ \ \text{ for } \frac14< r \le \frac12
    \vspace{5pt}     \\
     \varphi_t+(2s-1)\tilde{c}_t, \ \ \ \text{ for } \frac12<s\le 1,         
  \end{array} \right.
  \qquad
  \eeq
  with $\tilde{c}_t \gg 1$ sufficiently large to ensure that
  \beq
  \Ical_t(\varphi_t+\tilde{c}_t) < \Ical_t(v_{1,t}) -10.
  \eeq
  Clearly $\Gamma_t \in \Pcal$. Furthermore, by virtue of the above estimates \eqref{5-1}, \eqref{5-2}, \eqref{5-3}, for $t>0$ sufficiently small, we have:
  \bear
  \max_{s\in [0,1]} \Ical_t(\Gamma_t(s)) &\le& 2\max\limits_{c\ge 0}\{-e^ct^2\int_SKe^{\varphi_t}+e^{-c}\int_Se^{-\varphi_t}+4\pi(g-1)c \} \\
  &&\ \ \ \ \ \ \ \ +\frac12\|\nabla\varphi_t\|^2 + 4\pi(g-1)\int_S\varphi_t\\
  &\le& 8\pi\log\frac{1}{t^2} + C,
  \eear
  for some suitable $C>0$, independent of $t$. In view of \eqref{mtp1}, this proves \eqref{ap2a}.

  To obtain the reverse inequality, we decompose:
  \beq
  v_{2,t} = w_t+\ttc_t, \ \ \ \text{ with }\ttc_t = \fint_Sv_{2,t}(z) dA.
  \eeq
We use the Moser-Trudinger inequality (see for instance \cite{Aub82}) to estimate:
\bea\label{mti2}
 t^2\int_SKe^{w_t+\ttc_t} &\le&  t^2e^{\ttc_t}\|K\|_{\infty}\int_Se^{w_t} \notag \\
  &\le & t^2Ce^{\ttc_t}e^{\frac{\|\nabla w_t\|^2}{16\pi}}.
  \eea
By \eqref{5.1}, it is necessary that:
  \be\label{408a}
  \lim_{t\searrow 0}t^2\int_SKe^{v_{2,t}} \ge 4\pi,
  \ene
and so from \eqref{mti2} we find that, 
  \be\label{rev1}
  \|\nabla w_t\|^2 \ge 16\pi\log\frac{1}{t^2} - 16\pi\ttc_t - C_0,
  \ene
  for some suitable constant $C_0 > 0$.

  As a consequence, we find:
   \bea\label{rev2}
  8\pi\log\frac{1}{t^2} + 8\pi(g-2)\ttc_t - 2t^2\int_SKe^{v_{2,t}}+2\int_Se^{-v_{2,t}} -C_0&\le&\Ical_t(v_{2,t})\notag\\
  &\le&8\pi\log\frac{1}{t^2} +C.
  \eea
We are assuming $g \ge 3$, and also we know that: $\ttc_t>0$ and $t^2\int_SKe^{v_{2,t}}dA \in (0,4\pi(g-1))$. Thus, from 
  \eqref{rev2} we easily derive \eqref{average} and \eqref{ap2}.
\ep
%0----------------
Now we will prove Theorem ~\ref{ge3}:
\bp (of Theorem ~\ref{ge3})
Recall that we have set, $v_{2,t} = w_t+\ttc_t$, and from \eqref{average} and \eqref{ap2}, it follows that, as $t \searrow 0$,
\be\label{5.13a}
\frac{\|\nabla w_t\|^2}{|\log t^2|} \to 16\pi.
\ene
While by the first estimate in \eqref{mti2} and \eqref{408a}, we also have:
\be\label{5.13b}
\lim_{t\to 0}\frac{\log(\int_Se^{w_t})}{|\log t^2|} \ge 1, \text{ as } t\searrow 0.
\ene
As a consequence of \eqref{5.13a}, \eqref{5.13b} and the Moser-Trudinger inequality, as $t \searrow 0$, we find:
\be\label{4012}
\frac{\log(\int_Se^{w_t})}{\|\nabla w_t\|_2^2}\to \frac{1}{16\pi}.
\ene
Moreover, by \eqref{average} and Proposition ~\ref{unstable}, it follows that,
\be\label{4013}
\max_S w_t\to +\infty, \ \text{ as } t\to 0^{+}.
\ene
By the improved Moser-Trudinger inequality of Chen-Li \cite{CL91} (see Lemma 6.2.7 in \cite{Tar08} and also Malchiodi-Ruiz \cite{MR11}), 
and in view of \eqref{4012} and \eqref{4013}, we have that, there exists a unique point $p_0\in S$, such that, $\forall r>0$ sufficiently small, 
the following holds as $t \searrow 0$:
\be\label{4014}
\frac{\int_{B_r(p_0)} e^{w_t}}{\int_S e^{w_t}}\to 1,
\ene
\be\label{4015}
\max_{B_r(p_0)} w_t\to +\infty,
\ene
and
\be\label{4016}
\max_{S\backslash B_r(p_0)} w_t\le C_r,
\ene
for a suitable constant $C_r>0$. In other words, $p_0$ is the \underline{unique} blow-up point for $w_{t_n}$, along any sequence 
$t_n\searrow 0$, see \cite{BM91}. That is, if $p_n\in S$ satisfies:
\beq
w_{t_n}(p_n) = \max_Sw_{t_n} \to +\infty, \ \text{ as } n \to \infty
\eeq 
then:
\be\label{pt}
p_{n}\to p_0, \text{ as } n\to +\infty.
\ene
We shall show that, 
\be\label{4017}
K(p_0) \not= 0,
\ene
and since $K = \|\alpha\|_\sigma^2$, $p_0$ cannot be a zero for $\alpha \in Q(\sigma)$.

In order to see this, we use \eqref{average} and \eqref{408a} to find that,
\beq
\log(t^2\int_SKe^{w_t}) = -\ttc_t + O(1) = O(1), \text{ as } t\to 0.
\eeq
Thus, as $t \searrow 0$, we have:
\bear
O(1) &=&\Ical_t(v_{2,t})+8\pi\log\frac{1}{t^2} \\
&=& \frac12\|\nabla w_t\|^2-8\pi\log(\int_SKe^{w_t}) +8\pi(g-2)\log(t^2\int_SKe^{w_t}) +O(1),
\eear
that is, 
\be\label{4018}
\frac12\|\nabla w_t\|^2-8\pi\log(\int_SKe^{w_t}) =O(1).
\ene
On the other hand, from \eqref{4014} we easily check that,
\be\tag{5.22a}
\frac{\int_{S\backslash B_r(p_0)} Ke^{w_t}}{\int_{B_r(p_0)}e^{w_t}}\to 0, \text{ as } t\to 0,
\ene
and therefore, as $t \searrow 0$:
 \bear
 \log(\int_{S} Ke^{w_t}) &=&  \log(\int_{B_r(p_0)} Ke^{w_t}+ \int_{S\backslash B_r(p_0)} Ke^{w_t})\\
 &\le& \log\left(\max_{B_r(p_0)}(K)\right)+\log\int_{B_r(p_0)}e^{w_t}+\log\left(1+ \frac{\int_{S\backslash B_r(p_0)} Ke^{w_t}}{\int_{B_r(p_0)}e^{w_t}}\right) \\
  &<& \log\left(\max_{B_r(p_0)}(K)\right)+\log\int_{S}e^{w_t}+o(1).
 \eear
 As a consequence, from \eqref{4018} and the Moser-Trudinger inequality, as $t\to 0^{+}$, we find:
  \bear
C_1&\ge& \frac12\|\nabla w_t\|^2-8\pi\log(\int_SKe^{w_t}) \\
 &\ge& \frac12\|\nabla w_t\|^2-8\pi\log(\int_Se^{w_t})-8\pi\log\left(\max_{B_r(p_0)}(K)\right)+o(1) \\
  &>& -C_2-8\pi\log\max_{B_r(p_0)}(K) +o(1).
 \eear
 with suitable positive constants $C_1$ and $C_2$.

 Thus, we obtain:
 \beq
 \max_{z\in B_r(p_0)}(K(z)) \ge e^{-C}, \ \ \forall r>0,
 \eeq
 with a suitable constant $C>0$, independent of $r>0$. So by letting $r \searrow 0$, we get that $K(p_0) > 0$ 
 and \eqref{4017} is proved.

 At this point, for any sequence $t_n \searrow 0$, we can apply Theorem D for the sequence $v_{2,t_n}$. In view of \eqref{4015} and 
 \eqref{4016}, we know that $v_{2,t_n}$ can admit exactly \underline{one} blow-up point at $p_0$ with $K(p_0)\not=0$. Therefore \eqref{m1} must 
 hold with $m=1<g-1$, and consequently properties (i)-(iii) must hold for $v_{2,t_n}$.

 Since this holds along any sequence $t_n \searrow 0$, we obtain the desired conclusion.
 \ep
%%%%%%%%%%%%%%%%%%%%%%%%
\subsection{Blow-up analysis when $g=2$} 
When the surface is of genus $g =2$, the asymptotic behavior of $v_{2,t}$ is governed by the extremal properties of the Moser-Trudinger 
functional $\Jcal$ in \eqref{J} (see for instance \cite{Tru67, Mos71, Aub82}). Indeed, the goal of this subsection is to prove the following: 
\bt\label{g2bu}
Let $S$ be of genus $g=2$. Then, as $t\searrow 0$, we have:
\be\label{51}
t^2\int_SK(z)e^{v_{2,t}}dA\to 4\pi,
\ene  
\be\label{52}
\ttc_t =\fint_S v_{2,t}dA\to +\infty,
\ene  
\be\label{53}
\Ical_t(v_{2,t}) -8\pi\log\frac{1}{t^2} \to\inf_{w\in E}\Jcal(w) - 8\pi,
\ene  
where $\Jcal$ is the Moser-Trudinger functional defined in \eqref{J}, and $E = \{w\in H^1(S): \int_Sw(z)\ dA =0\}$.

Furthermore, by setting 
\be\label{w_t}
w_t = v_{2,t} - \ttc_t \in E,
\ene
then the following alternative holds:
\ben[{\bf (i)}]
\item \underline{either}, $\Jcal$ attains its infimum on $E$, and {\aas} $t = t_n\to 0$, as $n\to \infty$, we have:
\be\label{55a}
w_n \to w_0, \ \ \ \text{ uniformly in }\ C^{2,\beta}(S),
\ene
and 
\be\label{55b}
t^2e^{\ttc_t} \to \frac{4\pi}{\int_S Ke^{w_0}},
\ene
with $w_0$ satisfying the following equation on $S$:
\be\label{55c} \left\{
  \begin{array}{cl}
    -\Delta w_0 = 8\pi \left(\frac{K(z)e^{w_0}}{\int_S K(z)e^{w_0} dA}-\frac{1}{4\pi} \right) 
    \vspace{5pt} \\
    \Jcal(w_0) = \inf\limits_{w\in E}\Jcal(w).                
  \end{array} \right.
  \qquad
\ene
\item 
\underline{or}, the functional $\Jcal$ does not attain its infimum on $E$, and {\aas} $t = t_n\to 0$, as $n\to \infty$, for 
\beq
p_n\in S \ \text{ with } w_{t_n}(p_n) = \max_S w_{t_n},
\eeq
we have:
\be\label{5-6a}
 p_n \to p_0\in S, \ w_{t_n}(p_n)\to+\infty,
\ene
and 
\be\label{5-6}
t_n^2K(z)e^{v_{2,t_n}} \to 4\pi\delta_{p_0},
\ene
weakly in the sense of measure, and 
\be\label{5-7}
w_{t_n} \to 4\pi G(\cdot, p_0),
\ene 
uniformly in $C_{loc}^{2,\beta}(S\backslash\{p_0\})$, where $0<\beta<1$, with the blow-up point $p_0\in S$ satisfying:
\be\label{5-8}
4\pi\gamma(p_0,p_0)+\log K(p_0) =\max_{p\in S}\{4\pi\gamma(p,p)+\log K(p) \},
\ene
and in particular, $\alpha(p_0) \not= 0$.
\een
\et
\br
Clearly, if we knew the uniqueness of the minimum of the Moser-Trudinger functional $\Jcal$ on $E$ (when attained), or of the 
maximum point of the function $4\pi\gamma(p,p)+\log K(p)$, we could claim the convergence above as $t \to 0$, not only {\aas} 
$t=t_n\to 0$ as $n\to\infty$.
\er
Concerning the existence of a global minimum of $\Jcal$ in $E$, we briefly recall the work of Ding-Jost-Li-Wang (\cite{DJLW97}) 
and Nolasco-Tarantello (\cite{NT98}):
\bl\label{5.3} %For all $w \in H^1(S)$ with $\int_Sw(z)dA =0$, 
We have: 
\be\label{inf-J}
\inf_{w\in E}\Jcal(w) \le -8\pi(\max\limits_{p\in S}\{4\pi\gamma(p,p)+\log(K(p)\} + \log(2\pi(g-1))+1),
\ene
and the infimum is attained if \eqref{inf-J} holds with a strict inequality.
\el
\bp
See \cite{DJLW97, NT98}.
\ep
On the basis of Lemma ~\ref{5.3}, the existence of a global minimum for the extremal problem:
\be\label{DJLW}
\inf_{w\in E}\{\frac12\int_S|\nabla w|^2\ dA-8\pi\log(\fint_SK(z)e^wdA)\},
\ene
was ensured by the authors in \cite{DJLW97, NT98} under the following sufficient condition:
\be\label{suff}
\Delta_{g_\sigma}\log K(p_0) > -\left(\frac{8\pi}{|S|_\sigma}-2\kappa(p_0)\right)
\ene
with $p_0$ satisfying \eqref{5-8}, and $\kappa$ the Gauss curvature of $(S,\sigma)$. See also \cite{Tar08}.

For our geometrical problem, we have $K = \frac{|\alpha|^2}{det(g_\sigma)}$, with $g_\sigma$ the {\hym}, and $\alpha\in Q(\sigma)$ a 
{\hqd} on the {\RS} $(S,\sigma)$. So for any $p_0\in S$ with $\alpha(p_0)\not= 0$, we have $\kappa(p_0) = -1$, $|S|_\sigma =4\pi$ 
(note that $g=2$), and 
\beq
\Delta_{g_\sigma}\log K(p_0) = -4.
\eeq
Therefore we see that both sides of \eqref{DJLW} are equal to $-4$, and in this sense we just ``missed" to satisfy this sufficient condition 
\eqref{suff}.
\bp
We first apply \eqref{elelem} to the solution $v_{2,t}$ which (for $g=2$) implies that \eqref{m1} must hold with $m=g-1=1$. Therefore 
we have, as $t\searrow 0$,
\beq
t^2\int_SKe^{v_{2,t}} \to 4\pi, \ \ \text{ and } \int_Se^{-v_{2,t}} \to 0.
\eeq 
As before, we write $v_{2,t} = w_t + \ttc_t$, and (by Jensen's inequality) we find:
\beq
\ttc_t\to +\infty \ \ \text{as } t\searrow 0.
\eeq
This establishes \eqref{51} and \eqref{52}. Notice that we are now in the situation described by part (ii) of Theorem D.

In order to establish \eqref{53}, we use \eqref{51} and $g=2$, to conclude that the mean value $\ttc_t$ of $v_{2,t}$ must satisfy \eqref{hat-v} 
with the ``plus" sign. In other words, 
\be\label{5-11}
e^{\ttc_t} = \frac{2\pi+\sqrt{4\pi^2-(t^2\int_S Ke^{w_t} dA)(\int_Se^{-w_t}dA)}}{t^2\int_S Ke^{w_t} dA}.
\ene
So, we can use \eqref{5-11} to write 
\bea\label{5-12}
\Ical_t(v_{2,t}) &=&\frac12\|\nabla w_t\|_2^2 + 8\pi\log\frac{1}{t^2} -8\pi\log\fint_SKe^{w_t} +8\pi\notag \\
&&+8\pi\log(\frac{2\pi+\sqrt{4\pi^2-(t^2\int_S Ke^{w_t})(\int_Se^{-w_t})}}{4\pi}) -4t^2\int_SKe^{w_t+\ttc_t}  \notag\\
&=&\frac12\|\nabla w_t\|_2^2 + 8\pi\log\frac{1}{t^2} -8\pi\log\fint_SKe^{w_t} \notag\\
&&\ \ +8\pi\log(\frac{2\pi+\sqrt{4\pi^2-(t^2\int_S Ke^{w_t})(\int_Se^{-w_t})}}{4\pi}) \notag\\
&& \ \ \ \ -4(2\pi+\sqrt{4\pi^2-(t^2\int_S Ke^{w_t})(\int_Se^{-w_t})} ) + 8\pi.
\eea
Consequently,
\be\label{5.38a}
\Ical_t(v_{2,t}) - 8\pi\log\frac{1}{t^2} =\frac12\|\nabla w_t\|_2^2 -8\pi\log\fint_SKe^{w_t} -8\pi +o(1), \ \ \text{as } t \to 0,
\ene
and we derive the lower bound:
\be\label{5-12a}
\lim_{t\to 0}\left(\Ical_t(v_{2,t}) -8\pi\log\frac{1}{t^2}\right) \geq \inf_{w\in E}\Jcal(w) - 8\pi.
\ene  
To obtain the reversed inequality, we will construct some ``optimal" path. To this purpose, for any fixed $w\in E$, we find $t_w>0$ 
sufficiently small, such that: 
\beq
(t^2\int_S Ke^{w_t})(\int_Se^{-w_t})< 4\pi^2, \ \ \ \forall \ t \in (0, t_w).
\eeq 
So for every $t\in (0,t_w)$, we can define
\be\label{5-13}
\ttc_t^{\pm}(w) =\log\left(\frac{2\pi\pm\sqrt{(2\pi)^2-t^2\int_S K(z)e^w dA\int_Se^{-w}dA}}{t^2\int_S K(z)e^w dA}\right).
\ene
Also set, corresponding to the stable solution $v_{1,t}$: 
\be\label{5-14}
\ttc_{1,t} = \fint_S v_{1,t} \to 0, \ \ \text{as } t \searrow 0,
\ene 
and 
\be\label{5-15}
w_{1,t} = v_{1,t} -\ttc_{1,t} \to 0, \ \ \text{strongly in } C^{2,\beta}(S), \ \text{as } t \searrow 0.
\ene 
We define the following path:
   \be\label{path-w}
  \Gamma_{t,w}(s) =\left\{
  \begin{array}{cl}
    v_{1.t} -4sw_{1,t}, \ \ \ \text{ for } 0\le s \le \frac14
    \vspace{5pt} \\
   (4s-1)(w+\ttc_t^{-}(w)) +2(1-2s)\ttc_{1,t}, \ \ \ \text{ for } \frac14< s \le \frac12
    \vspace{5pt}     \\
     w+\ttc_t^{-}(w)+(2s-1)\tilde{C}_t, \ \ \ \text{ for } \frac12<s\le 1,         
  \end{array} \right.
  \qquad
  \ene
with $\tilde{C}_t > 0$ fixed sufficiently large (depending on $w$), to ensure that,
\beq
\Ical_t(w+\ttc_t^{-}(w)+\tilde{C}_t) < \Ical_t(v_{1,t}) -10, \ \ \forall \ t\in (0,t_w).
\eeq
Therefore $\Gamma_{t,w} \in \Pcal$, the path space defined in \eqref{pathspace}. Since 
\beq
\ttc_w^{-} \to \log\fint_Se^{-w}, \ \ \text{as } t \searrow 0,
\eeq
we readily check that,
\be\label{5-16}
\Ical_t(\Gamma_{t,w}(s)) \le C(w), \ \ \text{for }s \in [0, \frac12], \ \ t\in (0, t_w),
\ene 
with a suitable constant $C(w) > 0$ depending on $w$ only.

On the other hand, for $s \in [\frac12, 1]$ and $t\in (0, t_w)$, we have:
\bear
\Ical_t(\Gamma_{t,w}(s)) &\le&\frac12\|\nabla w_t\|_2^2 + 2\max_{c\ge \ttc_w^{-}}\{-t^2e^c\int_SKe^w+e^{-c}\int_Se^{-w}+4\pi c \} \\
&=&\frac12\|\nabla w_t\|_2^2 -2t^2e^{\ttc_t^{+}(w)}\int_SKe^w+ 2e^{-\ttc_t^{+}(w)}\int_Se^{-w} + 8\pi \ttc_t^{+}(w).
\eear
So, by observing that $v = w+\ttc_t^{+}(w)$ satisfies the integral identity \eqref{gb}, we can use \eqref{5-13} to show that 
\bea\label{5-17}
\max_{s\in [\frac12,1]}\Ical_t(\Gamma_{t,w}(s)) &\le&\frac12\|\nabla w\|_2^2 + 8\pi\log\frac{1}{t^2}  -8\pi\log\fint_SKe^{w} \notag\\
&&\ \ +8\pi\log\left(\frac{2\pi+\sqrt{4\pi^2-(t^2\int_S Ke^{w})(\int_Se^{-w})}}{4\pi}\right) +8\pi \notag\\
&& \ \ \ \ -4\left(2\pi+\sqrt{4\pi^2-(t^2\int_S Ke^{w})(\int_Se^{-w})}\right). 
\eea
So from \eqref{5-16} and \eqref{5-17}, for $t>0$ sufficiently small, we find:
\bea\label{5-18}
\Ical_t(v_{2,t}) -8\pi\log\frac{1}{t^2}&\le& \max_{s\in [0,1]}\Ical_t(\Gamma_{t,w}(s))  -8\pi\log\frac{1}{t^2} \notag\\
&\le& \frac12\|\nabla w\|_2^2 -8\pi\log\fint_SKe^{w} +8\pi \notag\\
&&\ \ +8\pi\log\left(\frac{2\pi+\sqrt{4\pi^2-(t^2\int_S Ke^{w})(\int_Se^{-w})}}{4\pi}\right) \notag\\
&& \ \ \ \ -4\left(2\pi+\sqrt{4\pi^2-(t^2\int_S Ke^{w})(\int_Se^{-w})}\right). 
\eea
As a consequence, we get:
\be\label{5-19}
\overline{\lim\limits_{t \searrow 0}}\left(\Ical_t(v_{2,t}) -8\pi\log\frac{1}{t^2}\right) \le \frac12\|\nabla w\|_2^2 -8\pi\log\fint_SKe^{w} -8\pi. 
\ene
Since \eqref{5-19} holds for every $w\in E$, and using \eqref{5-12a}, we establish \eqref{53}. Actually, from \eqref{5.38a} and 
\eqref{5-19}, we see that,
\be\label{5-20}
\lim\limits_{t \searrow 0}\left(\frac12\|\nabla w_t\|_2^2 -8\pi\log\fint_SKe^{w_t}\right) = \inf_E\Jcal. 
\ene
Next we wish to show that $w_t$ satisfies the ``compactness" alternative in part (ii) of Theorem D ({\aas} $t=t_n\searrow 0$) {\ifif} $\Jcal$ 
attains its infimum in $E$.

To this purpose, we fix $w\in E$, and as before set $t_w>0$ sufficiently small to ensure that, 
\beq
A_t(w) = (t^2\int_S Ke^{w_t})(\int_Se^{-w_t})< 4\pi^2, \ \ \ \forall \ t \in (0, t_w).
\eeq 
We set a function
\be\label{5-21}
f(A) = \log\left(\frac{2\pi+\sqrt{4\pi^2-A}}{4\pi}\right) + 8\pi -4(2\pi+\sqrt{4\pi^2-A}),
\ene
for $A \in [0,4\pi^2]$. Clearly this is a monotone increasing function of $A$ in $[0,4\pi^2]$.

From \eqref{5-18} and \eqref{5-12}, it follows that, for $\forall w\in E$, and $\forall t\in (0,t_w)$, there holds:
\be\label{5-22}
\frac12\|\nabla w_t\|_2^2 -8\pi\log\fint_SKe^{w_t} + f(A_t(w_t))\le \frac12\|\nabla w\|_2^2 -8\pi\log\fint_SKe^{w}+ f(A_t(w)).
\ene
Therefore if we assume the functional $\Jcal$ attains its infimum at $w_0$, that is, 
\be\label{5-22a}
\Jcal(w_0) = \inf_E\Jcal,
\ene
then we can use $w=w_0$ in \eqref{5-22} to conclude that, 
\be\label{5-23}
 f(A_t(w_t)) \le  f(A_t(w_0)), \ \ \forall \ t \in (0,t_{w_0}).
\ene
Now we use the fact that $f$ is increasing in $A$ and from \eqref{5-23} to derive that, 
\beq
(\fint_SKe^{w_t})(\fint_Se^{-w_t}) \le (\fint_SKe^{w_0})(\fint_Se^{-w_0}),
\eeq
with $\fint_Se^{-w_t} \ge 1$, by Jensen's inequality.

Therefore, for suitable $C_1>0$, we have,   
\be\label{5-24}
\fint_SKe^{w_t} \le  C_1, \ \ \forall \ t \in (0,t_{w_0}).
\ene
Using $v_{2,t} =w_t+\ttc_t$ in Lemma ~\ref{2-4}, we find that,
\be\label{5-25}
\fint_SKe^{w_t} \ge C_2,
\ene
for suitable $C_2>0$, and moreover $w_t$ is {\ub} in $W^{1,q}(S)$, $1< q<2$. Now from  \eqref{5-11}, \eqref{5-24} and \eqref{5-25}, 
we can deduce that:
\be\label{5-26}
\frac{1}{C_3} \le t^2e^{\ttc_t} \le C_3, \ \ \forall \ t \in (0,t_{w_0}),
\ene
with a suitable constant $C_3>1$.

In addition, from \eqref{5-20}, with $w=w_0$, we get
\bea\label{5-27}
\frac12\|\nabla w_t\|_2^2 &\le& 8\pi\log\fint_SKe^{w_t}+\inf_E\Jcal+ f(A_t(w_t))- f(A_t(w_0)) \notag \\
&\le&C_4,\ \ \forall \ t \in (0,t_{w_0}),
\eea
with a suitable constant $C_4>0$.

So in case the functional $\Jcal$ attains its infimum in $E$, then we can use the estimates in \eqref{5-26} and \eqref{5-27} together with 
elliptic estimates and well known regularity theory, to conclude that, $w_t$ is {\ub} in $C^{2,\beta}(S)$-norm, with $0 < \beta < 1$, 
for any $t \in (0, t_{w_0})$. Consequently, {\aas} $t_n\searrow 0$, $w_n:=w_{t_n}$ satisfies the ``compactness" property of part (ii) in 
Theorem D. In other words, \eqref{55a} holds with $w_0$ satisfying \eqref{55b} and \eqref{55c}.

Next suppose the functional $\Jcal$ does not attain its infimum in $E$. Therefore, $w_t$ can not satisfy the ``compactness" property in part 
(ii) in Theorem D. Consequently, by \eqref{51} we know that, along a sequence $t_n\searrow 0$, the sequence $w_n:=w_{t_n}$ must admit 
one $(m=1)$ blow-up point $p_0\in S$, satisfying \eqref{5-6a}, \eqref{5-6}, and \eqref{5-7}. So we are left to show that \eqref{5-8} holds. To 
this purpose, from \eqref{5-20}, we know that $w_n$ defines a minimizing sequence for $\Jcal$ in $E$ and $\max\limits_Sw_n\to+\infty$, i.e., 
blow-up occurs. Therefore, we can use for $w_n$, the estimates detailed in \cite{DJLW97} and \cite{NT98} for any blow-up minimizing 
sequences of $\Jcal$, to show that, 
\bea\label{5-28}
\inf_E\Jcal &=& \lim_{n\to +\infty}\frac12\|\nabla w_n\|_2^2-8\pi\log\fint_SKe^{w_n} \notag \\
&\ge& -8\pi\left(4\pi\gamma(p_0,p_0)+\log K(p_0) + \log\frac{\pi}{|S|} +1\right).
\eea
On the other hand, when $\Jcal$ does not attain its infimum, we also know that, 
\be\label{5-29}
\inf_E\Jcal = -8\pi\max_{p\in S}\left(4\pi\gamma(p,p)+\log K(p) + \log\pi +1\right).
\ene
See Lemma ~\ref{5.3} and \cite{DJLW97}. Now \eqref{5-8} follows immediately from \eqref{5-28} and \eqref{5-29}.
\ep
%%%%%%%%%%%%%%%%%%%%%%%%%%%%%%%%%%%%%%%%
\section{Prescribing extrinsic curvature}
In this section, we wish to investigate the possibility to obtain a {\mi} of $S$ into a {\htm} with prescribed \underline{total extrinsic curvature}. 
Namely, for given $\rho\in (0, 4\pi(g-1))$, we require that for the induced metric $g_0$ we have:
\be\label{6-1}
\rho = \int_S (det_{g_0}\Pi)dA_{g_0}. %\in (0, 4\pi(g-1)).
\ene
%%%%%%%%%%%%%%%%%%%
\subsection{Main result and three lemmata}
Our main result is the following:
\begin{thx}\label{main5}
Fixing a {\cs} $\sigma \in T_g(S)$, and a {\hqd} $\alpha \in Q(\sigma)$, and $\rho\in (0, 4\pi(g-1))\backslash\{4\pi m, m=2, \cdots, g-2\}$, there 
exists a constant $t_\rho \in (0, \tau_0]$ ($\tau_0 = \tau_0(\sigma,\alpha)> 0$ given in Theorem ~\ref{Uh}), such that $S$ admits a {\mi} of 
data $(\sigma, t_\rho\alpha)$ into some {\htm}, with corresponding total extrinsic curvature satisfying \eqref{6-1}. 
\end{thx}
In order to establish this result, we need to provide a solution $v_\rho$ for the problem $(1)_{t_\rho}$, for some $t_\rho \in (0, \tau_0]$ satisfying:
\be\label{6-2}
t_\rho^2\int_SK(z)e^{v_\rho}\ dA = \rho, \ \ \ K =\|\alpha\|_\sigma^2 = \frac{|\alpha|^2}{det(g_\sigma)}.
\ene
To this purpose, we recall from \S3.1 the Mean Field formulation of the problem $(1)_t$, as described in Lemma ~\ref{mfe}. Then for given 
$\rho \in (0, 4\pi(g-1))$ we need to find a solution $w$ of the equation \eqref{w}, that is:
 \beq \left\{
  \begin{array}{cl}
    -\Delta w = 2\rho\left(\frac{K(z)e^{w}}{\int_S K(z)e^{w} dA}-\frac{1}{|S|} \right) +2(4\pi(g-1)-\rho)\left(\frac{e^{-w}}{\int_Se^{-w} dA}-\frac{1}{|S|} \right)
    \vspace{5pt} \\
    \int_Sw(z)dA =0.                    
  \end{array} \right.
  \qquad
\eeq
We call this equation, (namely \eqref{w}), the problem $(3)_\rho$.

To this end, we start with the following observation:
\bl\label{lemma6-1}
If $w$ solves the problem $(3)_\rho$ with $\rho\in (0, 4\pi(g-1))$, then:
\be\label{6-7}
\|\frac{e^{-w}}{\int_Se^{-w}dA}-\frac{1}{|S|} \|_{L^\infty} < \frac{1}{4\pi(g-1)-\rho}.
\ene
\el
\bp
We recall from Lemma ~\ref{mfe} that when $w$ solves problem $(3)_\rho$, then we define:
\be\label{6-5}
\ttc_\rho =\log\left(\frac{\int_Se^{-w} dA}{4\pi(g-1)-\rho} \right),
\ene
and 
\be\label{6-6}
t_\rho^2 =\frac{\rho(4\pi(g-1)-\rho)}{(\int_S K(z)e^{w} dA)(\int_Se^{-w} dA)},
\ene
and we see that $v_\rho = w+\ttc_\rho$ is a solution to problem $(1)_{t_\rho}$. Since we have:
\beq
\left(4\pi(g-1)-\rho\right)\left(\frac{e^{-w}}{\int_Se^{-w}dA}-\frac{1}{|S|}\right) = e^{-v_\rho} - \fint_Se^{-v_\rho}\ dA,
\eeq
and $v_\rho > 0$, then, $\forall \ z\in S$:
\beq 
|e^{-v_\rho(z)} - \fint_Se^{-v_\rho}\ dA| < 1-e^{-v_\rho}(z) < 1,
\eeq
and this establishes \eqref{6-7}.
\ep
As a consequence of Corollary ~\ref{3.5}, we know that, for any $\rho \in (4\pi(m-1), 4\pi m)$, $m=\{1, \cdots, g-1\}$, the 
Leray-Schauder degree of the Fredholm operator associated to the problem $(3)_\rho$ is well-defined, and its value only depends 
on $m$. To be more precise, we recall that the Laplace-Beltrami operator $\Delta = \Delta_{g_\sigma}$ is invertible on the space $E$. 
We denote by: 
\beq
(\Delta|_E)^{-1}: E \to E
\eeq
its (smooth) inverse. Thus each solution to the problem $(3)_\rho$ corresponds to a zero of the following operator:
\be\label{6-8}
F_\rho(w) = w + T_\rho^0(w) + B_\rho(w), \ \ \ \forall \ w\in E,
\ene
with 
\be\label{6-9}
T_\rho^0(w) =2\rho(\Delta|_E)^{-1}\left(\frac{Ke^{w}}{\int_SKe^{w}dA}-\frac{1}{|S|} \right),
\ene
and 
\be\label{6-10}
B_\rho(w) =2\left(4\pi(g-1)-\rho\right)(\Delta|_E)^{-1}\left(\frac{e^{-w}}{\int_Se^{-w}dA}-\frac{1}{|S|} \right).
\ene
Therefore, in view of Lemma ~\ref{lemma6-1}, there exists a suitable constant $C>0$ (independent of $\rho$), such that if 
$w \in E$ is a solution of problem $(3)_\rho$, that is $F_\rho(w) =0$, then $\|B_\rho(w\|\le C$.

As a consequence,  for any $\rho \in (4\pi(m-1), 4\pi m)$, with $m=\{1, \cdots, g-1\}$, we find a radius $R_\rho$ sufficiently large, such that, 
for each $t\in [0,1]$ it is well-defined at zero the Leray-Schauder degree $d_{\rho, t}$ of the operator
\be\label{F-t}
F_\rho^t(w) = w + T_\rho^0(w) + tB_\rho(w), 
\ene
in the ball $B_{R_\rho} = \{w\in E: \|w\|\le R_\rho\}$. Moreover, by the homotopy invariance of the Leray-Schauder degree, we have 
\beq
d_\rho = d_{\rho, t}, \ \ \forall\ t\in [0,1].
\eeq 
In particular, 
\be\label{6-11}
d_\rho = d_{\rho, 0},
\ene
where $d_{\rho, 0}$ is the Leray-Schauder degree of the operator 
\beq
F_\rho^0(w) = w + T_\rho^0(w), \ \ \ \forall \ w\in E,
\eeq
whose zeroes correspond to solutions of the following problem: 
\beq\left\{
  \begin{array}{cl}
    -\Delta w = 2\rho \left(\frac{K(z)e^{w}}{\int_S K(z)e^{w} dA}-\frac{1}{|S|}\right), \ \ \text{on } S
    \vspace{5pt} \\
    \int_Sw(z)dA =0,                    
  \end{array} \right.
  \qquad
\eeq
Actually for every $\rho\not\in 4\pi \mathbb{N}$, the Leray-Schauder degree of the operator $F_\rho^0$ has been computed by Chen-Lin 
(\cite{CL03, CL15}), exactly when the weight function $K$ admits isolated zeros (which is our case) $\{q_1, \cdots, q_N\}$, each 
with integer multiplicity $\nu(q_j), \ 1\le j \le N$. More precisely we have the following:
\bl\label{lemma6-3}(\cite{CL15})
If $\nu(q_j) \in \mathbb{N} \ \forall 1\le j \le N$, and the genus of the surface is greater than zero, then $d_{\rho, 0}> 0$.
\el
\bp
See Corollary 1.2 in \cite{CL15}.
\ep
%%%%%%%%%%%%%%%%%%%%%
\subsection{Proof of Theorem E}
We now complete the proof of Theorem E:
\bp
Since in our case, the weight function $K(z) = \|\alpha\|_\sigma^2$ with $\alpha \in Q(\sigma)$ a {\hqd} on $S$, we know that $\alpha$ 
admits isolated zeroes of integer multiplicity and total number (counting multiplicity) equal to $4(g-1)$. Thus, we can apply 
Lemma ~\ref{lemma6-3} together with \eqref{6-11} to conclude that, for every $\rho \in (0, 4\pi(g-1))\backslash\{4\pi m, m=1, \cdots, g-2\}$, 
the Leray-Schauder degree $d_{\rho}> 0$. In other words, for such range of $\rho$'s, we know that the problem $(3)_\rho$ admits at least 
one solution. To complete the proof, we need to show that, when $g\ge 3$, then we have the existence of a solution for problem $(3)_\rho$ 
also when $\rho = 4\pi$.

To this purpose, we once again exploit the work of Chen-Lin in \cite{CL10, CL02}. We take a sequence $w_n$ of the solutions 
to problem $(3)_{\rho_n}$, with $4\pi m\not =\rho_n\to 4\pi m$, for some $m\in\{ 1, \cdots, g-2\}$. We assume that, as $n\to+\infty$, the following holds:
\be\label{6-12}
w_n\rightharpoonup w_0 \ \ \ \text{ weakly in } W^{1,q}(S), 1< q<2,
\ene
and 
\be\label{6-13}
\max_S w_n \to +\infty,
\ene
Indeed, in case $w_n$ was {\ub} in $S$, then by elliptic estimates ({\aas}) it would converge to a solution to problem $(3)_{\rho=4\pi m}$, and 
for $m=1$ we would obtain our solution in this way. Thus, we assume \eqref{6-13} and we want to establish a sign for the quantity $\rho_n-4\pi m$. 
This delicate task has been carried out by Chen-Lin in \cite{CL10, CL02} for the sequence $z_n = w_n-\zeta_n$, satisfying \eqref{6*} and \eqref{4.16a}, 
with $\zeta_n$ defined in \eqref{zeta} and satisfying ({\aas}):
\be\label{6-14}
\zeta_n \to \zeta_0, \ \ \ \text{strongly in }C^{2,\beta}(S), \ \text{as } n \to +\infty,
\ene
with $\zeta_0$ the unique solution for:
\be\label{6-15} \left\{
  \begin{array}{cl}
    -\Delta \zeta_0 = 8\pi(g-m-1)\left(\frac{e^{-w_0}}{\int_Se^{-w_0}}-\frac{1}{|S|}\right)\ \ \ \text{ on }S,
    \vspace{5pt} \\
    \int_S\zeta_0(z)dA =0,               
  \end{array} \right.
  \qquad
\ene
From \eqref{6-13} we know that, $\max\limits_S z_n \to +\infty$. Therefore, by using Theorem ~\ref{star}, $z_n$ must admit a finite number of 
blow-up points, say $\{p_1, \cdots, p_s\}\subset S$, for which \eqref{4.17c} holds with $m=\sum\limits_{j=1}^s(1+n(p_j))$ and 
\beq
w_0(z) = \zeta_0(z) + 8\pi\sum\limits_{j=1}^s(1+n(p_j))G(z,p_j).
\eeq

If we further assume that these blow-up points are not zeroes of the weight function $K =\|\alpha\|_\sigma^2$, that is 
$\alpha(p_j) \not= 0, \ \forall j=1, \cdots, s$; then $n(p_j) = 0$ and $m =s$. In this situation, Chen-Lin in \cite{CL15} were able to control the exact 
decay to zero of the quantity: $\rho_n-4\pi m$. In particular they showed that, the sign of $\rho_n-4\pi m$ is the same of the following quantity:
\be\label{6-16}
 \sum_{j=1}^md_j\left(\Delta\log h^{*}(p_j)+\frac{8\pi m}{|S|}-2\kappa(p_j) \right),
\ene
where $d_j$'s are suitable (positive) constants, $h^{*}=Ke^{\zeta_0}$, and $\kappa$ is the Gauss curvature of $S$. Take into account also that the 
expression \eqref{6-16} was given in \cite{CL15} by formulae $(2.3)$ and $(2.10)$, written under the normalization $|S| =1$. Now, for 
$p\in S: \alpha(p)\not= 0$, by means of \eqref{6-15}, we compute:
\bear
\Delta\log(Ke^{\zeta_0})(p)  + \frac{8\pi m}{|S|} -2\kappa(p) &=& \Delta\log\|\alpha\|_\sigma^2 + \Delta\zeta_0+\frac{2m}{g-1}+2\\
&=&-4-8\pi(g-m-1)\left(\frac{e^{-w_0}}{\int_Se^{-w_0}}-\frac{1}{4\pi(g-1)} \right)\\
&& \ \ \ \ +\frac{2m}{g-1} + 2\\
&=&-8\pi(g-m-1)\frac{e^{-w_0}}{\int_Se^{-w_0}} \\
&<& 0, \ \ \ \ \ \forall\ m=\{1,\cdots, g-2\}.
\eear
Therefore, we may conclude that, if $K$ (and hence $\alpha$) does not vanish at the blow-up points of $w_n$, then for $n$ sufficiently 
large, we have: $\rho_n - 4\pi m < 0$. That is, blow-up can only occur from the ``right".

Since for $m=1$, the solutions to problem $(3)_{\rho_n}$ with $\rho_n\to 4\pi$ can admit only \underline{one} blow-up point $p_0\in S$ 
which must satisfy $K(p_0)\not= 0$. Therefore, we can use the information above, to see that for $\rho_n>4\pi$ and $\rho_n\searrow 4\pi$, 
the corresponding solution $w_n$ cannot blow-up, and so ({\aas}) it converges to the desired solution of $(3)_{\rho=4\pi}$.
\ep
We conclude the section with two remarks.
\br 
When the sequence of solutions $w_n$ to problem $(3)_{\rho_n}$ with $\rho_n\to 4\pi m$, blows up at a zero of $K=\|\alpha\|_\sigma^2$, we 
suspect that similar information about the sign of the quantity $\rho_n - 4\pi m$ should hold. This is confirmed by the more involved analysis 
developed in \cite{CL10}, where the authors provide sharp estimates about the behavior of the sequence $z_n$ of \eqref{4.16a}, \eqref{6*}, 
which blows up at a \underline{zero} of the weight function $K$, but only when such a zero is of {\underline{non-integer}} multiplicity.
\er
\br
Finally we note that, in view of \eqref{1-c}, \eqref{1-d}, by choosing $\alpha \in Q(\sigma)$ with zeroes of multiplicity greater than $g-2$, we 
can always guarantee that blow-up never occurs at its zeroes.
\er

%%%%%%%%%%%%%%%%%%%%%%%%%%%%%%
\bibliographystyle{amsalpha}
\bibliography{ref-HLT}
\end{document}